\documentclass[11pt]{amsart}

\newtheorem{theorem}{Theorem} [section]
\newtheorem{prop}[theorem]{Proposition}
\newtheorem{lemma}[theorem]{Lemma}
\newtheorem{cor}[theorem]{Corollary}

\usepackage{amsfonts}
\usepackage{amsopn}
\usepackage{epsfig}
\usepackage{color}
\usepackage{amsmath}
\usepackage{graphicx}
\usepackage{amssymb}

\newcommand\C{{\bf C}}
\newcommand\Chat { {\hat{\C}} } 

\renewcommand\P{{\bf P}}
\newcommand\R{{\bf R}}
\newcommand\Z{{\bf Z}}
\newcommand\D{{\bf D}}
\newcommand\Hyp{{\bf H}}

\newcommand\del{\partial}

\newcommand\eps{\varepsilon}

\renewcommand\phi{\varphi}

\newcommand\iso{\simeq} 

\renewcommand\O{\mathcal{O}} 

\newcommand\M{\mathcal{M}}
\newcommand\scriptC{\mathcal{C}}

\newcommand\Aut{\operatorname{Aut}}
\newcommand\SL{\operatorname{SL}}
\newcommand\PSL{\operatorname{PSL}}

\renewcommand\H{\operatorname{H}}

\newcommand\graph{\operatorname{Graph}}
\renewcommand\mod{\operatorname{mod}}  
\newcommand\Id  {\operatorname{Id}}     

\renewcommand\gcd {\operatorname{gcd}} 
\newcommand\Res {\operatorname{Res}} 

\newcommand\Rat  {\operatorname{Rat}} 
\newcommand\Ratbar {\overline{\Rat}}  
\newcommand\Mbar {\overline{M}}

\begin{document}

\noindent

\title{The boundary of the moduli space of quadratic 
     rational maps}  
\author{Laura DeMarco}

\date{December 21, 2004}

\begin{abstract}
Let $M_2$ be the space of quadratic rational maps $f:\P^1\to\P^1$, 
modulo the action by conjugation of
the group of M\"obius transformations.   
In this paper a compactification $X$ of $M_2$ is defined, as a
modification of Milnor's
$\Mbar_2\iso\C\P^2$, by choosing representatives of 
a conjugacy class $[f]\in M_2$ such that the measure of maximal 
entropy of $f$ has conformal barycenter at the origin in $\R^3$, and
taking the closure in the space of probability measures.  
It is shown that $X$ is the smallest compactification of $M_2$ such 
that all iterate maps $[f]\mapsto [f^n]\in M_{2^n}$ extend continuously
to $X \to \Mbar_{2^n}$, where $\Mbar_d$ is the natural compactification
of $M_d$ coming from geometric invariant theory.  
\end{abstract}

\maketitle
\thispagestyle{empty}

 \section{Introduction}

For each $d\geq 2$, let $M_d = \Rat_d/\PSL_2\C$ denote the
space of degree $d$ rational maps $f:\Chat\to\Chat$, modulo 
the action by conjugation of the group of M\"obius transformations.
The moduli space is a complex orbifold of dimension $2d-2$.  
Iteration defines a sequence of regular maps 
$$\Phi_n: M_d \to M_{d^n},$$
given by $[f]\mapsto [f^n]$, where $[f]$ denotes the conjugacy
class of $f\in\Rat_d$.  

The aim of this paper is 
to define a compactification of the moduli space which is natural
from the point of view of dynamics, and in particular, one on 
which the iterate maps are well-defined.  
Two approaches to this end are presented here, one using results 
in geometric invariant theory and one
in terms of measures of maximal entropy.  In degree $d=2$, 
the two approaches are shown to be equivalent.   

\bigskip\noindent{\bf A formal solution.}  
In \cite{Silverman}, Silverman
studied a compactification $\Mbar_d$ of $M_d$, for each $d\geq 2$, 
by computing the stability criteria for the 
conjugation action of $\SL_2\C$ 
on $\Rat_d\hookrightarrow \P^{2d+1}$, according to Mumford's 
geometric invariant theory.  The iterate maps $\Phi_n$, however, 
do not define regular maps from  $\Mbar_d$ to $\Mbar_{d^n}$ 
for any $d\geq 2$ and $n\geq 2$ (see \S\ref{higher}).  

It is possible to define a compactification 
of the moduli space $M_d$ on which iteration is well-defined, by resolving 
the indeterminacy of each
rational iterate map $\Mbar_d \dashrightarrow \Mbar_{d^n}$ 
and passing to an inverse limit.  Namely, we can define $\Gamma_n$ to be the
closure of $M_d$ as it sits inside the finite product $\Mbar_d \times \Mbar_{d^2} \times
\cdots \times \Mbar_{d^n}$ via the first $n$ iterate maps $(\Id, \Phi_2, \ldots , \Phi_n)$.
There is a natural projection from $\Gamma_{n+1}$ 
to $\Gamma_n$ for every $n$, so we may 
take the inverse limit over $n$, 
  $$\hat{M}_d = \lim_{\longleftarrow} \Gamma_n.$$
The moduli space $M_d$ is a dense open subset of $\hat{M}_d$, where a conjugacy class
$[f]\in M_d$ is identified with the sequence $([f], [f^2], [f^3], \ldots )$
in $\hat{M}_d$.  The iterate
map $\Phi_n: M_d \to M_{d^n}$ extends continuously to $\hat{M}_d \to \hat{M}_{d^n}$,
by sending the sequence $([f], [f^2], [f^3], \ldots )$ to the sequence $([f^n], [f^{2n}], [f^{3n}],
\ldots )$.  It remains to understand the structure of this space and if there exists a 
concrete model for $\hat{M}_d$.

\bigskip\noindent{\bf Maximal measures and the barycenter.}
Given a rational map $f\in\Rat_d$, let $\mu_f$ denote the unique probability 
measure on $\Chat$ of maximal entropy 
\cite{Lyubich:entropy},\cite{FLM},\cite{Mane:unique}.
The support of $\mu_f$ is equal to the Julia set of $f$, and the measure
is invariant under iteration, 
$\mu_{f^n} = \mu_f$ for all $n\geq 1$.
The conformal barycenter of $\mu_f$ is its
hyperbolic center of mass, where the 
the unit ball in $\R^3$ is chosen as a model for $\Hyp^3$, and the unit
sphere $S^2$ is identified with the Riemann sphere $\Chat$ via stereographic projection
\cite{Douady:Earle} (see \S\ref{barycenter}).  

For each conjugacy class $[f]\in M_d$, 
we can choose a barycentered representative $f\in\Rat_d$,
one such that the conformal barycenter of $\mu_f$  is at the 
origin in $\R^3$.  The representative is unique up to the action of the compact
group of rotations $SO(3) \subset \PSL_2\C$.  
If $BCM$ denotes
the space of barycentered probability measures on $\Chat$ (with the
weak-$*$ topology), then we obtain a continuous map, 
  $$BC: M_d \to BCM/SO(3).$$  
Let $\overline{BCM}$ denote the closure of $BCM$ 
in the space of all probability measures, and consider the closure of the graph of $BC$, 
 $$X_d = \overline{\graph(BC)} \subset \Mbar_d \times \overline{BCM}/SO(3).$$ 
This defines a new compactification of the moduli space $M_d$. 

\bigskip\noindent{\bf Quadratic rational maps.} 
In degree $d=2$, Milnor showed that the moduli space, $M_2 = \Rat_2/\PSL_2\C$, 
is an orbifold with underlying complex manifold 
isomorphic to $\C^2$ \cite[Lemma 3.1]{Milnor:quad}.  The
compactification $\Mbar_2 \iso \P^2$ has a boundary consisting of 
the conjugacy classes of degree 1 rational maps and one degree 0 map.
It is isomorphic to the geometric invariant theory compactification in degree 2
\cite[Thm 1.5]{Silverman}.

The main theorem of this paper shows that  
the compactification $X_2$ of $M_2$
by barycentered measures is an explicit model for the
formal construction of $\hat{M}_2$ which resolves the iterate maps.   

\begin{theorem}  \label{homeo}
The compactifications $\hat{M}_2$ and $X_2 = \overline{\graph(BC)}$
of $M_2$ are canonically homeomorphic. 
\end{theorem}

\noindent
In other words, there is a homeomorphism $\hat{M}_2\to X_2$ which restricts
to the identity on the moduli space $M_2$.   It is not true in general 
that the compactifications $\hat{M}_d$ and $X_d$ are homeomorphic for 
every $d\geq 2$.
Examples are given in Section \ref{higher}.  
However, for $d=2$ we have the following corollary. 

\begin{cor}  
The iterate maps $\Phi_n: M_2\to M_{2^n}$ extend continuously
to $X_2 \to X_{2^n}$ for every $n\geq 1$.
\end{cor}

\proof
Let $\{[f_k]\}_{k=0}^\infty$ be a sequence in $M_2$ such that
$[f_k] \to p \in \del M_2 \subset X_2$ as 
$k\to \infty$.  By the definition of $X_2$, there exist representatives
$f_k\in\Rat_2$ with barycentered measures of maximal entropy
$\mu_{f_k}$ which converge weakly to a probability measure 
$\nu$ as $k\to \infty$.
By Theorem \ref{homeo}, the point $p$ is identified with a unique
point in $\del M_2\subset \hat{M}_2$, and therefore it has well defined 
iterates $p^n\in \Mbar_{2^n}$ for all $n\geq 1$.  By the continuity  
of $\Phi_n : \hat{M}_2 \to \Mbar_{2^n}$ and the iterate-invariance of the 
measures, $\mu_{f_k^n} = \mu_{f_k}$, the sequence of iterates
$\Phi_n([f_k])$ must converge in $X_{2^n}$ to the point $(p^n, \nu)$.
\qed

\bigskip\noindent{\bf The structure of $X_2 = \hat{M}_2$.}  
As an inverse limit construction, the space $\hat{M}_2$ could have 
very complicated structure.  In fact, the boundary of $M_2$ in this space
can be fully understood.  The first result says $\hat{M}_2$ can not be 
embedded into any finite dimensional projective space.  

\begin{theorem}  \label{infinite}
No finite sequence of blow-ups of $\Mbar_2 \iso \P^2$ is enough to 
resolve all of the rational iterate maps $\Phi_n: \Mbar_2 \dashrightarrow
\Mbar_{2^n}$ simultaneously.
\end{theorem}

On the other hand, the space $\Gamma_n$, which is the closure of the 
graph of $(\Phi_2$, $\Phi_3$,  $\ldots, \Phi_n)$ in 
$\Mbar_2\times\Mbar_4\times \cdots\times \Mbar_{2^n}$ 
has a fairly simple structure, described
completely in Section \ref{blowup}.  In particular, there are no ideal points in 
the inverse limit space $\hat{M}_2 = \lim \Gamma_n$:  

\begin{theorem}  \label{finite}
Every sequence in $\hat{M}_2 \subset \prod_1^\infty \Mbar_{2^n}$ is 
determined by finitely many entries.
\end{theorem}  

\begin{figure}[htbp]
\begin{center}
\input{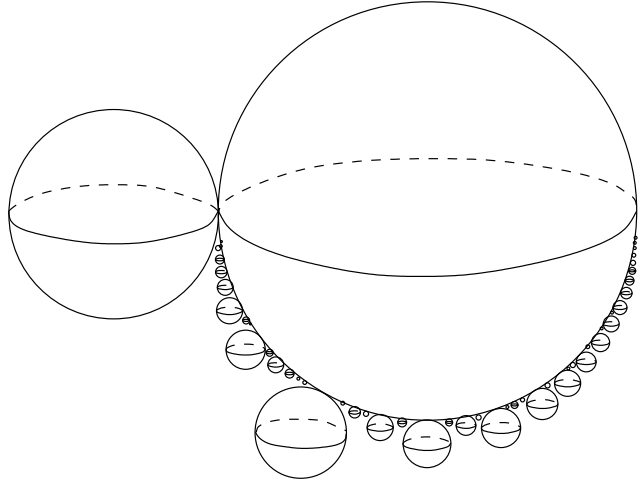}
\end{center}
\caption{Boundary of $M_2$ in $X_2 = \hat{M}_2$.}
\label{figure1}
\end{figure}

\noindent
Topologically, the boundary of $M_2$ in $\hat{M}_2$ is obtained 
from the boundary of $M_2$ in $\Mbar_2 \iso\P^2$ by successively attaching 2-spheres
at a countable collection of points in $\Mbar_2$, as in Figure \ref{figure1}.  

\bigskip\noindent{\bf Outline of the paper.}
The notation is fixed in Section \ref{ratd}, and a summary of results
from \cite{D:measures} is provided.  
The geometric invariant theory compactification 
$\Mbar_d$ is defined in Section \ref{GITstability}.  
In Section \ref{proper}, we show
that the iterate map $M_d\to M_{d^n}$ is proper.  Sections \ref{degree2}
and \ref{degree2proofs} are devoted to the study of iteration in 
degree 2, and we give the proof of Theorem \ref{infinite}.  Section 
\ref{blowup} contains a study of the structure of $\hat{M}_2$ and the
proof of Theorem \ref{finite}.  The space of barycentered measures 
$BCM/SO(3)$ is studied in Section \ref{barycenter}, and the proof 
of Theorem \ref{homeo} is contained in Section \ref{maintheorem}.
Section \ref{higher} is devoted to a study of the iterate 
map $M_d\to M_{d^n}$ in general degrees $d\geq 2$.  
Some concluding remarks about the definitions of $\hat{M}_d$
and $X_d$ are given in Section \ref{final}.  

\bigskip\noindent{\bf Acknowledgements.}  
The analysis of quadratic rational maps used here appeared
first in \cite{Milnor:quad} and, in greater detail, in \cite{Epstein:bounded}
where Epstein studied the structure of hyperbolic components in 
$M_2$ and gave the first examples of discontinuity of the iterate map 
at the boundary.  In fact, as shown in Sections \ref{degree2} and \ref{degree2proofs},
Epstein's examples are the only examples which demonstrate discontinuity
in degree 2.  
The geometric invariant theory approach relies on the 
results in \cite{Silverman}.  I am grateful to A.~Epstein, J.~Harris, 
J.~Hubbard, C.~McMullen, and J.~Milnor 
for helping me formulate the results in this paper.

\bigskip
\section{$\Ratbar_d$ and the probability measures at the boundary} 
\label{ratd}

In this section, we fix notation and terminology.  
We state some facts about the iterate map $\Rat_d \to 
\Rat_{d^n}$ and the measures of maximal entropy
 from \cite{D:measures}.  

\bigskip\noindent{\bf The compactification $\Ratbar_d$.}
Let $\Rat_d$ denote the space of holomorphic maps $f:\Chat\to\Chat$
of degree $d$ with the topology of uniform convergence.
For each $d$, $\Rat_d$ is naturally identified with 
the complement of a hypersurface in 
  $$\Ratbar_d = \P\H^0(\P^1\times\P^1, \O(d,1)),$$
by sending 
$f\in\Rat_d$ to the section which vanishes along the graph of $f$.
The space of rational maps $\Rat_d$ is therefore a smooth, affine
variety, and 
we obtain an isomorphism 
$\Ratbar_d \iso \P^{2d+1}$.  Alternatively, 
each point $f\in\Rat_d$ determines
a pair of degree $d$ homogeneous polynomials, unique up to scale,
  $$f(z:w) = (P(z,w):Q(z,w)),$$ 
and the space of such pairs 
is $\P^{2d+1}$, parametrized by the coefficients of $(P,Q)$.
In particular, 
  $$\Rat_d \iso \P^{2d+1} - V(\Res),$$
where $V(\Res) = \{(P,Q): \Res(P,Q) = 0\}$ is the resultant hypersurface.
In $\Ratbar_d$, the hypersurface $V(\Res)$ corresponds to the collection
of all sections with reducible zero locus.  

Given a pair $(P,Q)$, 
the zeroes of the homogeneous polynomial $H = \gcd(P,Q)$, as points in $\P^1$, 
will be called the {\bf holes} of  the associated $f\in\Ratbar_d$ and the 
multiplicity of
a zero the {\bf depth} of the hole.  Each $f\in\Ratbar_d$ determines
a holomorphic map 
$$\phi_f = (P/H: Q/H): \P^1\to \P^1$$ 
of degree $\leq d$.  
We will often write 
  $$f = (P:Q) =  H_f\phi_f$$ 
where $H_f = \gcd(P,Q)$.

\bigskip\noindent{\bf Coordinates on $\P^1$.}
A point $(z:w)\in\P^1$ will regularly be identified with 
$z/w\in\hat{\C}$.  Any distances on $\Chat$ will be measured
in the spherical metric.  A ball of radius $r$ about a point 
$p\in\Chat$ will be denoted $B(p,r)$.  

\bigskip\noindent{\bf The measure of maximal entropy.}
Fix $d\geq 2$.
Given a rational map $f\in \Rat_d$, there is a unique (non-exceptional)
probability measure $\mu_f$ on $\Chat$  such that 
  $$\frac{1}{d} f^* \mu_f = \mu_f$$
\cite{Lyubich:entropy},\cite{FLM},\cite{Mane:unique}.  The
measure $\mu_f$ is of maximal entropy ($\log d$) with support
equal to the Julia set of $f$.  Ma\~n\'e showed that 
the function $f\mapsto \mu_f$ is continuous from $\Rat_d$
to the space of probability measures with the weak-$*$ topology
\cite{Mane}.  

\bigskip\noindent{\bf Iteration on $\Ratbar_d$.}
We provide here a summary of relevant definitions and 
statements from \cite{D:measures}.  
The main object of study in \cite{D:measures} is the relation
between the iterate maps, $f\mapsto f^n$, extended to $\Ratbar_d$
and the extension of the map of maximal measures, $f \mapsto \mu_f$,
to $\Ratbar_d$.  

The {\bf indeterminacy locus} $I(d) \subset \Ratbar_d$ is the collection of 
$f = H_f\phi_f$ with $\deg \phi_f = 0$ and such that the constant 
value of $\phi_f$ is a hole of $f$.  It has codimension $d+1$
in $\Ratbar_d$.  See Figure \ref{figure2}.

\bigskip
\begin{figure}[htbp]
\begin{center}
\input{Ratbar.pstex_t}
\end{center}
\caption{Graphs in $\P^1\times\P^1$ of $f\in \Ratbar_3$:  (a) $f\in\Rat_3$, 
(b) $f = H_f\phi_f\in\del\Rat_3$ with $\deg\phi_f = 2$, (c) $f\in I(3)$.}
\label{figure2}
\end{figure}
\bigskip

\begin{theorem}  \label{Ratindet}
\cite[Thm 2]{D:measures}
The indeterminacy locus of the iterate map $\Ratbar_d 
\dashrightarrow \Ratbar_{d^n}$ is $I(d)$ for all $n\geq 2$.  
\end{theorem}

\noindent
Consequently, any element $f = H_f\phi_f \in\Ratbar_d - I(d)$
has well-defined forward iterates $f^n$ for all $n\geq 2$.  A 
direct computation yields the formula, 
  $$f^n = \left( \prod_{k=0}^{n-1} (\phi^{k*}H_f)^{d^{n-k-1}} \right) \phi_f^n$$
\cite[Lemma 7]{D:measures}.  

\bigskip\noindent{\bf Atomic probablity measures.}
For each $f = H_f\phi_f \in \del\Rat_d$ such that 
$\deg \phi_f > 0$, 
a purely atomic
probability measure $\mu_f$ is defined by the following triple sum,
  $$\mu_f = \sum_{n=0}^\infty \frac{1}{d^{n+1}} 
\sum_{\{H_f(h) = 0\}}  \sum_{\{\phi_f^n(z) = h\}}\delta_z,$$
where the middle sum is over all holes of $f$, the inner sum is 
over all preimages of those holes, and the outer sum is over 
all iterates of $\phi_f$, all counted with mulitiplicity.  Because
the number of holes is $d - \deg\phi_f$, it is easy to check that $\mu_f$ has 
total mass one.  For $\deg\phi_f = 0$, we define the probability measure by 
  $$\mu_f = \frac{1}{d} \sum_{\{H_f(h) = 0\}} \delta_h.$$
One can define pull-back of measures by 
any $f\not\in I(d)$, and the measure $\mu_f$ is the unique probability
measure satisfying 
$f^*\mu_f = d\cdot\mu_f$ \cite[Prop 10]{D:measures}.  

\begin{theorem}  \label{limitmu}
\cite[Thm 1(a)]{D:measures}
Given any sequence $\{f_k\}$ in $\Rat_d$ converging
to $f \in \del\Rat_d - I(d)$ in $\Ratbar_d$,
the measures of maximal entropy $\mu_{f_k}$ 
converge weakly to $\mu_f$.  
\end{theorem}

We will use the following three lemmas throughout this text.
The first two follow directly from a comparison of the formula
for an iterate of $f\in\del\Rat_d$ with the definition of $\mu_f$.
For $f = H_f\phi_f\in\Ratbar_d$, let $d_h(f)$ denote 
the depth of $h\in\P^1$ as a hole of $f$ and let $m_h(\phi_f)$
be the multiplicity of $z=h$ as a solution to $\phi_f(z) = \phi_f(h)$.
Note that $m_h(\phi_f) = 1$ if and only if $h$ is not a critical 
point of $\phi_f$.  By convention, $m_h(\phi) = 0$ for all $h$ if 
$\phi$ is constant, and the 0-th iterate $\phi^0$ is the identity
map.  

\begin{lemma} \label{Lemma5} \cite[Lemma 5]{D:measures}
For each $f = H_f\phi_f \in \Ratbar_d$ and $z\in\Chat$, we have 
  $$\mu_f(\{z\}) = \frac{1}{d} \sum_{n=0}^\infty 
\frac{m_z(\phi_f^n) d_{\phi_f^n(z)}(f)}{d^n}.$$
\end{lemma}

\begin{lemma} \label{Cor8} \cite[Cor 8]{D:measures}
For each $f\in\del\Rat_d - I(d)$, the depths of the holes of the iterates 
of $f$ are given by 
 $$d_z(f^n) = d^{n-1} \cdot d_z(f) + \sum_{k=1}^{n-1} d^{n-1-k} 
m_z(\phi_f^k) d_{\phi_f^k(z)}(f).$$
Therefore, the sequence $\{d_z(f^n)/d^n: n\geq 1\}$
is non-decreasing, and 
  $$\mu_f(\{z\}) = \lim_{n\to\infty} \frac{d_z(f^n)}{d^n}.$$
\end{lemma}

\begin{lemma} \label{sequence} \cite[Lemmas 14, 15]{D:measures}
Suppose $\{f_k\}$ is a sequence in $\Rat_d$ 
converging to $f = H_f\phi_f = (P:Q)$ in $\Ratbar_d$.  
\begin{itemize}
\item[(i)]  The sequence of rational maps $f_k$ converges to $\phi_f$ 
locally uniformly on the 
complement of the holes of $f$ in $\Chat$, and 
\item[(ii)] if $f$ has a hole at $h$ of depth $d_h$ and neither $P$ 
nor $Q$ is $\equiv 0$, then any neighborhood of $h$ contains 
at least $d_h$ zeroes and poles of $f_k$ (counted with multiplicity)
for all sufficiently large $k$.  
\end{itemize}
\end{lemma}

\medskip
We also need some more general results on the structure
of the composition map.  Recall the notation from Lemma 
\ref{Lemma5}.

\begin{lemma} \label{composition}
The composition map 
 $$\scriptC_{d,e}: \Ratbar_d\times\Ratbar_e\dashrightarrow \Ratbar_{de},$$
which sends a pair $(f,g)$ to the composition $f\circ g$, is 
continuous away from 
  $$I(d,e) = \{(f,g)=(H_f\phi_f,H_g\phi_g): \phi_g\equiv c
\mbox{ and } H_f(c)=0 \}.$$
Furthermore, for each 
$(f,g) \in \Ratbar_d\times\Ratbar_e$ such that $\deg\phi_g>0$,  
 $$d_z(f\circ g) = d\cdot d_z(g) + m_z(\phi_g)\cdot d_{\phi_g(z)}(f).$$
\end{lemma}

\proof
In the coordinates on $\Ratbar_d$ and $\Ratbar_e$ 
given by the coefficients of $f$ and $g$, the
composition map is defined by polynomial functions, so 
it suffices to show that 
$\scriptC_{d,e}(f,g)$ is well-defined for each pair $(f,g)\not\in I(d,e)$.
Write $f = (H_fP_f:H_fQ_f) = H_f\phi_f\in\Ratbar_d$
and $g = (H_gP_g:H_gP_g)  = H_g\phi_g \in\Ratbar_e$.  Let $d' = \deg \phi_f$.
The composition $f\circ g$ can be computed directly by 
\begin{eqnarray*}  
 \scriptC_{d,e}(f,g) &=& (H_f (H_gP_g,H_gQ_g) P_f(H_gP_g, H_gQ_g) : \\
& &   \qquad  \qquad H_f(H_gP_g, H_gQ_g) Q_f (Q_gP_g, H_gQ_g) ) \\
&=& (H_g)^{d-d'}H_f(P_g,Q_g) (H_g)^{d'} 
(P_f(P_g,Q_g): Q_f(P_g,Q_g))  \\
&=& (H_g)^d H_f(P_g,Q_g) \: \phi_f\circ \phi_g.
\end{eqnarray*}
By hypothesis, neither $H_f(P_g,Q_g)$ nor $(H_g)^d$ vanishes identically
and $\phi_f\circ\phi_g$ is a well-defined rational map.  

The formula for the depth of the holes of the composition
$f\circ g$ follows directly from the formula for 
$\scriptC_{d,e}(f,g)$ above.  
\qed

\bigskip
\section{The GIT stability conditions}  \label{GITstability}

The action of $\SL_2\C$ by conjugation on $\Rat_d$ extends to
$\Ratbar_d = \P\H^0(\P^1\times\P^1, O(d,1))$ by 
the diagonal action on $\P^1\times\P^1$.  In this section
we describe the stability conditions for this action according
to geometric invariant theory (GIT), computed in 
\cite{Silverman}.  (See also \cite{Mumford:GIT}.)  We relate this notion of stability
to the atomic probability measures $\mu_f$ for $f\in\del\Rat_d$,
defined in \S\ref{ratd}.

Silverman  showed that the moduli space 
$M_d = \Rat_d/\PSL_2\C$ 
exists as a geometric quotient scheme which is 
affine, integral, connected, and of finite type over $\Z$.  
Furthermore, $\Mbar_d$ is a geometric quotient for $d$ 
even and a categorical quotient for $d$ odd, and it is 
proper over $\Z$ \cite[Thm 2.1]{Silverman}.  His computations
led to the following proposition:  

\begin{prop} \label{Silvermanprop} \cite[Prop 2.2]{Silverman}
A point $f\in\Ratbar_d$ is stable (respectively, semistable) for
the conjugation action of $\SL_2\C$ if and only if 
there are no elements in the conjugacy class of $f$ of the form 
  $$(a_0 z^d + a_1 z^{d-1} w + \cdots + a_d w^d: 
b_0 z^d + b_1 z^{d-1} w + \cdots + b_d w^d)$$
with $a_i = 0$ for all $i < (d-1)/2$ (respectively, $i\leq (d-1)/2$) 
and $b_j=0$ for all $j < (d+1)/2$ (respectively, $j\leq (d+1)/2$).  
\end{prop}

Denote the set of stable points by $\Rat^s_d \subset\Ratbar_d$
and the semistable points by $\Rat^{ss}_d$, and note that
$\Rat_d^s = \Rat_d^{ss}$ if and only if $d$ is even.  
Therefore, the compact GIT quotients are defined by 
$\Mbar_d = \Rat_d^s/\PSL_2\C$ for $d$ even
and $\Mbar_d = \Rat_d^{ss}//\PSL_2\C$ for $d$ odd. 
Roughly speaking, $\Rat_d^s$ is the largest open
$\PSL_2\C$-invariant subset of $\Ratbar_d$ in which 
all $\PSL_2\C$-orbits are closed, 
and so the quotient space $\Rat_d^s/\PSL_2\C$ is Hausdorff.   
When $d$ is odd, elements of $\Rat^{ss}_d$ represent
the same point in $\Mbar_d$
if the closures of their orbits intersect in $\Rat_d^{ss}$.

\bigskip\noindent{\bf The stable and semistable points.}
The following is a reformulation of Proposition \ref{Silvermanprop}
in the language of this paper.  Let $f = H_f\phi_f$ and $g = 
H_g\phi_g$ be two elements in $\Ratbar_d$.  Then $f$ and $g$ 
are in the same $\PSL_2\C$-orbit if and only if there exists 
$A\in\PSL_2\C$ such that $\phi_g = A\phi_fA^{-1}$ and the holes
of $g$ are the image under $A$ of the holes of $f$ (and of corresponding
depths).  
For even $d\geq 2$,  a point $f=H_f\phi_f\in\Ratbar_d$ is stable
(or semistable) if 
\begin{itemize}
\item[(i)] the depth of each hole is $\leq d/2$, and
\item[(ii)] if the depth of $h\in\P^1$ is $d/2$ then 
$\phi_f(h) \not= h$.
\end{itemize}
For odd $d\geq 3$, a point $f\in\Ratbar_d$ is stable if
\begin{itemize}
\item[(i)] the depth of each hole is $\leq (d-1)/2$, and
\item[(ii)] if the depth of $h\in\P^1$ is $(d-1)/2$ then 
$\phi_f(h) \not= h$,
\end{itemize}
and $f\in\Ratbar_d$ is semistable if
\begin{itemize}
\item[(i)] the depth of each hole is $\leq (d+1)/2$, and
\item[(ii)] if the depth of $h\in\P^1$ is $(d+1)/2$ then 
$\phi_f(h) \not= h$.
\end{itemize}

\bigskip\noindent {\bf Instability of stability.}
The property of stability
is not generally preserved by iteration.  For example,
consider the point $h = (zw:z^2)\in\Ratbar_2$ which is stable
($h$ has one hole at $z=0$ and $\phi_h(z) = 1/z$). 
This point
$h$ does not lie in the indeterminacy locus
$I(2)$ and therefore has a well-defined
second iterate, namely, 
  $$h^2 = (z^3w:z^2w^2)\in\Ratbar_4.$$
The second iterate $h^2$ coincides with the identity map away 
from a hole 
of depth 2 at $z=0$ and a hole of depth 1 at $z=\infty$.  
Consequently, $h^2$ is {\em not} stable.  

As another example, note that any degenerate polynomial
is either unstable itself or will eventually be unstable
after iteration.  That is, for any $0<k< d$, consider
  $$p = (w^k Q(z,w): w^d)$$
where $Q$ is a homogeneous polynomial of degree $d-k$
such that $Q(1,0)\not=0$,
so that $\phi_p(z)$ is the polynomial $Q(z,1)$ and $p$
has a hole of depth $k$ at $z=\infty$.  The iterates of 
$p$ are of the form, 
  $$p^n = (w^{kd^{n-1} + k(d-k)d^{n-2} + \cdots + k(d-k)^{n-1}} 
        Q^n(z,w): w^{d^n}),$$
and $p^n$ is unstable for all $n$ such that $\deg Q^n = (d-k)^n$ is 
less than or equal to $d^n/2$.  

For the examples just mentioned, we can compute the associated
probability measures using Lemma \ref{Lemma5},
 $$\mu_h = \frac{2}{3} \delta_0 + \frac{1}{3} \delta_\infty$$
and $\mu_p = \delta_\infty$.  The following propositions show 
that we can read from the measures that some iterate
of $h$ and $p$ will be unstable.

\begin{prop}  \label{evendeg}
Suppose $d$ is even and $f\not\in I(d)$.  Then 
$f^n \in \Ratbar_{d^n}$ is stable for all $n\geq 1$ if and only if
$\mu_f(\{z\}) \leq 1/2$ for all $z\in\P^1$.
\end{prop}

\proof
Let $d_z(f^n)$ denote the depth of $z$ as a hole of $f^n$. 
Write $f=H_f\phi_f$ and note that $\phi_{f^n} = \phi_f^n$ (see
\S\ref{ratd}).  
From Lemma \ref{Cor8}, the hypothesis on $\mu_f$ 
implies that $d_z(f^n) \leq d^n/2$ for all $z$ and all $n\geq 1$.
Suppose for some $n$ and 
$z$ we have $d_z(f^n) = d^n/2$  and $\phi_f^n(z) = z$.  The 
depth of $z$ as a hole of the composition  $f^{2n} = f^n\circ f^n$ 
must satisfy  
$d_z(f^{2n}) \geq  d^n(d^n/2) + 1 > d^{2n}/2$ by Lemma \ref{composition}, 
since $z$ is also 
one of the preimages of the hole at $z$ for $f^n$, providing
a contradiction.  

Conversely, suppose that $f^n$ is stable for all $n$.  Then, 
in particular, $d_z(f^n) \leq d^n/2$ for all $z$ and so again by 
Lemma
\ref{Cor8}, $\mu_f(\{z\}) = \lim_{n\to\infty} d_z(f^n)/d^n$ 
cannot exceed $1/2$. 
\qed  

\begin{prop} \label{odddeg}
Suppose $d$ is odd and $f\not\in I(d)$.
Then  $f^n\in\Ratbar_{d^n}$ is semi-stable for all $n\geq 1$ if and only if 
$\mu_f(\{z\}) \leq 1/2$ for all $z\in\P^1$. 
Furthermore,  if $\mu_f(\{z\}) < 1/2$ 
for all $z\in\P^1$, then $f^n$ is stable for all $n\geq 1$.  
\end{prop}

\proof  Let $d_z(f^n)$ denote the depth of $z$ as a hole of $f^n$,
and write $f= H_f\phi_f$.  
By Lemma \ref{Cor8}, $\mu_f(\{z\}) \leq 1/2 $
implies that $d_z(f^n) \leq (d^n -1)/2$ for all $n$ and $z$.  This
gives the second statement immediately: $f^n$ is semi-stable
for all $n$.  Conversely, if $f^n$ is semi-stable for all $n$,
then $d_z(f^n) \leq (d^n + 1)/2$ for all $n$.  But in the limit,
this implies that $\mu_f(\{z\}) \leq 1/2$.  

To prove the final statement, note again that $\mu_f(\{z\}) < 1/2$
implies that $d_z(f^n) \leq (d^n -1)/2$ for all $n$.
Suppose that for some $n$,  $d_z(f^n) = (d^n-1)/2$ and
$\phi_f^n(z) = z$.  Then, by Lemma \ref{Lemma5} applied to $f^n$, 
  $$\mu_{f^n}(\{z\}) \geq \frac{1}{d^n} \sum_{l=0}^\infty
\frac{(d^n-1)/2}{d^{nl}} = \frac{1}{2},$$  
which is a contradiction since $\mu_{f^n} = \mu_f$.
\qed

\bigskip
Note that the converse to the second statement of Proposition
\ref{odddeg} is false.  There exist odd degrees $d$ and 
$f\in\Ratbar_d$ for 
which $f^n$ is stable for all $n\geq 1$ but with $\mu_f(\{z\}) = 1/2$
for some point $z\in\P^1$.  Consider, for example, 
$f = (z^4w: zw^4) \in\Ratbar_5$
which has holes at $0$ and $\infty$, each of depth 1, and $\phi_f(z)
= z^3$.  Using Lemmas \ref{Lemma5} and \ref{Cor8}, it
is straightforward to compute that $\mu_f(0) = \mu_f(\infty) = 1/2$
while the depth of each of the two holes for $f^n$ is 
$(5^n - 3^n)/2 < (5^n -1)/2$.

\bigskip
\section{Properness of the iterate map on $M_d$} 
\label{proper}

The iterate map $\Rat_d \to \Rat_{d^n}$, given by $f\mapsto f^n$,
is a regular map between smooth, affine varieties.
It is $\PSL_2\C$-equivariant since $(\phi f \phi^{-1})^n = 
\phi f^n \phi^{-1}$ for all M\"obius transformations $\phi$
and all $n\geq 1$.  Therefore it descends to a regular map
on the affine moduli spaces, 
  $$\Phi_n: M_d \to M_{d^n}.$$
Iteration on $\Rat_d$ is proper (the preimage of 
any compact set in $\Rat_{d^n}$ is compact) if and only if 
the degree $d$ is at least 2
\cite[Cor 3]{D:measures}.  This implies the following.

\begin{prop}
The iterate map $\Phi_n: M_d \to M_{d^n}$ is proper
for every $d\geq 2$ and $n\geq 1$.  
\end{prop}

\proof
If $\{[f_k]\}_{k=0}^\infty$ is an unbounded sequence 
in $M_d$, then every sequence of representatives $f_k\in
\Rat_d$ is unbounded.  By properness of iteration on $\Rat_d$, 
every sequence 
of iterates $\{(f_k)^n\}_{k=0}^\infty$ is unbounded
in $\Rat_{d^n}$ \cite[Cor 3]{D:measures}.  
Consequently, $\Phi_n([f_k]) = [(f_k)^n]$ is 
unbounded in $M_{d^n}$.  
\qed

\bigskip
Since $M_d$ is an open dense subset of the GIT compactification
$\Mbar_d$, the map $\Phi_n$ defines a rational map on the 
closures, 
  $$\Phi_n : \Mbar_d \dashrightarrow \Mbar_{d^n}.$$
In Section \ref{higher}, we prove that 
$\Phi_n$ is not regular on $\Mbar_d$ for 
any $d\geq 2$ and $n\geq 2$ (Theorem \ref{discontinuity}).  
We aim to study the indeterminacy of this rational map 
in the boundary of $M_d$.  The following lemma begins
to address the relationship between the indeterminacy of 
$\Phi_n$ on $\Mbar_d$, the stability conditions in $\Ratbar_d$,
and the indeterminacy locus $I(d)\subset \del\Rat_d$ for the
general degree $d\geq 2$.

\begin{lemma} \label{stableindet}
Suppose $f\in\Ratbar_d$ satisfies $f\not\in I(d)$ and 
$f^n$ is stable for some $n>1$.  Then 
\begin{itemize}
\item[(i)]  $f$ is stable, and 
\item[(ii)] the iterate map $\Phi_n$ is continuous at $[f]\in\Mbar_d$.
\end{itemize}
\end{lemma}


\proof
Write $f = H_f\phi_f$.  
For $z\in\P^1$, let $d_z(f^n)$ denote the depth of $z$ 
as a hole of $f^n$.  The sequence 
$\{d_z(f^n)/d^n: n\geq 1\}$ is non-decreasing by 
Lemma \ref{Cor8}.

Suppose first that the degree $d$ is even.  
Stability of $f^n$ implies that $d_z(f^n) \leq d^n/2$ for all $z$, 
and therefore $d_z(f^k) \leq d^k/2$ for all $k\leq n$ and all $z$. 
If $d_z(f^n) = d^n/2$, then $\phi_f^n(z) \not= z$, and therefore 
$\phi_f(z)\not=z$
so we see that $f$ is stable.  It therefore determines a well-defined
point $[f]$ in $\Mbar_d$.  

If $d$ is odd, stability of $f^n$ implies that $d_z(f^n) \leq (d^n-1)/2$ for all
$z$.  Therefore $d_z(f)/d \leq (d^n-1)/2d^n < 1/2$ and so 
$d_z(f) \leq (d-1)/2$ for all $z$.  Furthermore, if we have
$d_z(f) = (d-1)/2$ with $\phi_f(z) = z$, then 
  $$d_z(f^n) \geq d^{n-1} \sum_0^{n-1} \frac{d-1}{2d^k}
        =  \frac{d^n - 1}{2}$$
and $\phi_f^n(z) = z$, which contradicts stability.  

Recalling that stability in $\Ratbar_d$ is an open condition,
we see that if $[f_t]$ is a family in $\Mbar_d$ converging
to $[f]$, then there is a family of stable representatives 
$f_t$ converging to $f$ in $\Ratbar_d$.  Since $f\not\in I(d)$,
the iterates $(f_t)^n$ converge to $f^n$ in 
$\Ratbar_{d^n}$.  Also, $(f_t)^n$ must be stable 
for all sufficiently small $t$, and therefore the iterates 
of $[f_t]$ converge in $\Mbar_{d^n}$ to $[f^n]$.  Therefore, 
$\Phi_n$ is continuous at $[f]\in\Mbar_d$.
\qed

\bigskip\noindent{\bf Question.}  Is the converse true, 
in the sense that if $f\in I(d)$ is stable, then $\Phi_n$ is 
indeterminate at $[f]$?  
And if $f\not\in I(d)$ is stable but $f^n$ is not stable, 
then is $\Phi_n$ indeterminate at $[f]$?

\bigskip
The answer is yes in degree $d=2$, as we shall 
see in the following sections.  What makes degree 2 particularly 
easy for computation is the following observation.  

\begin{lemma}  
The intersection $\Rat^{ss}_d \cap ~I(d)$ in $\Ratbar_d$
is empty if and only if $d = 2$.  
\end{lemma}

\proof
The statement is immediate from the definition of $I(d)$ and 
Proposition \ref{Silvermanprop}.   
\qed

\bigskip
\section{The moduli space $M_2\iso\C^2$}
\label{degree2}

In this section, we collect some fundamental facts about the 
moduli space of quadratic rational maps, and we describe completely
the indeterminacy locus of the iterate map 
$$\Phi_n : \Mbar_2 \dashrightarrow \Mbar_{2^n}$$ 
which sends $[f]$ to $[f^n]$.  We give the proof of Theorem 
\ref{infinite}.  

The analysis in this section is based
on the work of Milnor and Epstein \cite{Milnor:quad}, \cite{Epstein:bounded}
and the isomorphism between Milnor's compactification $\Mbar_2\iso\P^2$ 
and the geometric invariant theory compactification of $M_2$ 
\cite[Thm 1.5]{Silverman}.

\bigskip\noindent{\bf The compactification $\Mbar_2\iso\P^2$.} 
Every rational map $f:\Chat\to\Chat$ of degree 2 has three fixed 
points (solutions to $f(z)=z$), counted with multiplicity.  
The derivative of $f$ evaluated at a 
fixed point is called the {\bf multiplier} of the
fixed point. 
The multipliers of the fixed points are the solutions to a 
unique monic polynomial,
  $$x^3 - \sigma_1 x^2 + \sigma_2 x - \sigma_3 = 0,$$
and the residue formula (applied to the form $dz/(z-f(z))$) 
implies the relation $\sigma_3 = \sigma_1 - 2$. 
As the multipliers are conjugacy invariant, the $\sigma_i$ define 
functions on the moduli space $M_2 = \Rat_d/\PSL_2\C$.   
Milnor showed that the pair $(\sigma_1, \sigma_2)$ 
naturally parametrizes the moduli space, 
defining an isomorphism $M_2 \iso \C^2$ \cite[Lemma 3.1]{Milnor:quad}.
Consequently, a sequence $\{[f_k]\}$ is unbounded in 
$M_2$ if and only if some fixed point multiplier of $f_k$ tends to infinity.  

The boundary of $M_2$ in the compactification   
$\Mbar_2 \iso \P^2$, arising naturally from Milnor's isomorphism, 
corresponds to unordered triples of fixed point multipliers of 
the form $\{a, 1/a, \infty\}$ for $a\in\Chat$.  These triples can
be identified with  
the conjugacy classes of degree 1 and constant maps 
of the form $z\mapsto az+1$.  Indeed, under the identification 
of Milnor's $\Mbar_2$ with the GIT compactification 
(\cite[Thm 1.5]{Silverman}),
the line at infinity is parametrized in a two-to-one fashion by  
$a\mapsto [\Lambda_a] = [\Lambda_{1/a}]$
where $[\Lambda_a]$ is the class of  the following points in $\Ratbar_2$:
\begin{equation} \label{lambda}
\Lambda_a(z:w) =  \left\{ \begin{array}{l}
(az(z-w):w(z-w)), \mbox{ for } a \in \hat{\C} - \{0,1,\infty\},   \\
((z+w)(z-w):w(z-w)), \mbox{ for } a =1,  \\
((z+w)(z-w): 0), \mbox{ for } a =\infty, \mbox{ and } \\
(0: (z+w)(z-w)), \mbox{ for } a=0.  
\end{array} \right. 
\end{equation}
That is, $\Lambda_1$ is the parabolic M\"obius transformation $z\mapsto z+1$ 
with a hole at $z=1$, $\Lambda_\infty$ is the constant infinity map with holes
at $1$ and $-1$, $\Lambda_0$ is the constant 0 with holes at 
$1$ and $-1$, and for each $a\not= 0,1,\infty$, $\Lambda_a$ is given by 
$z\mapsto az$ with hole at $z=1$.
Recall by Lemma \ref{sequence} that any sequence in $\Rat_2$ converging
to $\Lambda_a$ in $\Ratbar_2$ will converge to the corresponding degree 0 or 1
map, locally uniformly on the complement of the holes of $\Lambda_a$.  

In this section we prove, 


\begin{theorem}  \label{indet}
For each $n\geq 2$, iteration defines a rational map
  $$\Phi_n: \Mbar_2 \dashrightarrow \Mbar_{2^n}$$
with indeterminacy locus given by
  $$I(\Phi_n) = \{ [\Lambda_a]\in\del M_2: a\not= 1 \mbox{ and } 
         a^q = 1 \mbox{ for some } 1 < q \leq n \}.$$
In particular, we have $I(\Phi_2) \subset I(\Phi_3) \subset I(\Phi_4)
\cdots$.
\end{theorem}

\bigskip\noindent{\bf Proof of Theorem \ref{infinite}.}  From 
Theorem \ref{indet}, we know that the indeterminacy set 
of $\Phi_n:\Mbar_2\dashrightarrow \Mbar_{2^n}$ is strictly increasing
with $n$.  Therefore, no finite sequence of blow-ups over points
in $I(\Phi_n)$ will suffice to resolve the indeterminacy of all
iterate maps simultaneously.  
\qed

\bigskip  
For $a\in\Chat$, 
let $\Lambda_a\in\Ratbar_2$ be defined by (\ref{lambda}).  
Checking the stability conditions of \S\ref{GITstability}, we see
that each of the points $\Lambda_a$ is stable, and these
elements
represent all of the stable conjugacy classes in $\Ratbar_2-\Rat_2$.  
Note also that $\Lambda_a\not\in I(2)$ for every $a\in\Chat$ and therefore
all forward iterates $\Lambda_a^n\in\Ratbar_{2^n}$ are well-defined by Theorem \ref{Ratindet}.  

\begin{lemma} \label{stableLambda}
The iterates $\Lambda_a^n\in\Ratbar_{2^n}$
are stable for all $n\geq 1$ if and only if $a\in\Chat$ is not
a primitive $q$-th root of unity for any $q\geq 2$.   For any
primitive $q$-th root of unity $\zeta$, the iterate $\Lambda_\zeta^n$ 
is stable if and only if $n< q$.  
\end{lemma}  

\proof
By the definition of the atomic probability measures $\mu_{\Lambda_a}$
given in Section \ref{ratd},
we have 
  $$\mu_{\Lambda_0} = \mu_{\Lambda_\infty} = \frac{1}{2} 
             \delta_1 + \frac{1}{2} \delta_{-1},$$
  $$\mu_{\Lambda_1} = \frac12 \sum_{k=0}^\infty \frac{1}{2^k} \delta_{1-k},$$
and
  $$\mu_{\Lambda_a} = \frac12 \sum_{k=0}^\infty \frac{1}{2^k} \delta_{1/a^k}$$
for all $a\not= 0,1, \infty$.  
We see immediately that $\mu_{\Lambda_a}(\{z\}) \leq 1/2$ for all $z\in\P^1$
if only if $a$ is not a root of unity.  By Proposition \ref{evendeg}, 
all iterates of $\Lambda_a$ must then be stable.

Now let $\zeta$ be a primitive $q$-th root of unity for some $q>1$.  
Writing $\Lambda_\zeta = H_\zeta \phi_\zeta$, note that $\phi_\zeta^k(z)
= \zeta^k z$, so that $\phi_\zeta^q$ is the identity map.  For each $n<q$
and $l = 0,1,\ldots,q-1$, it is easy to compute from Lemma \ref{Cor8}
that the depths of the holes of $\Lambda_\zeta^n$ are given by
 $$d_{1/\zeta^l}(\Lambda_\zeta^n) = 2^{n-1-l} \leq  2^n/2,$$
and $1/\zeta^l$ is not fixed by $\phi_\zeta^n$, so that 
$\Lambda^n_\zeta\in\Ratbar_{2^n}$ is stable.  
On the other hand, we find that $d_1(\Lambda_\zeta^q) = 2^{q-1}$ with 
$\phi_\zeta^q(1) = 1$ so that $\Lambda_\zeta^q$ is not stable.  
For each $n>q$, $d_1(\Lambda_\zeta^n) > 2^{q-1}$ so that 
$\Lambda_\zeta^n$ is not stable.  
\qed

\bigskip\noindent{\bf Epstein's normal forms}.
We now follow \cite{Epstein:bounded}.
Suppose $f\in\Rat_2$ has distinct
fixed points at 0, $\infty$, and 1, with multipliers 
$\alpha$, $\beta$, and $\gamma = (2-\alpha-\beta)/
(1 - \alpha\beta)$, respectively.  Then $f$ can be written
\begin{equation} \label{Nf}
f_{\alpha, \beta}(z) = z \frac{(1-\alpha)z + \alpha(1-\beta)}
                        {\beta(1-\alpha)z + (1-\beta)}.
\end{equation}
If the two critical points
of $f$ are distinct from the fixed point at $1$, then
$f$ is conjugate to 
\begin{equation}  \label{NF}
  F_{\gamma, \delta}(z) = \frac{\gamma z}{z^2 + \delta z + 1}
\end{equation}
for some $\delta\in \C$, where $F$ has critical points
at $1$ and $-1$ and a fixed point of multiplier $\gamma$ at 0.   Interchanging
the labelling of the critical points replaces $\delta$ with $-\delta$.  

\begin{figure}[htbp]
\begin{center}
\input{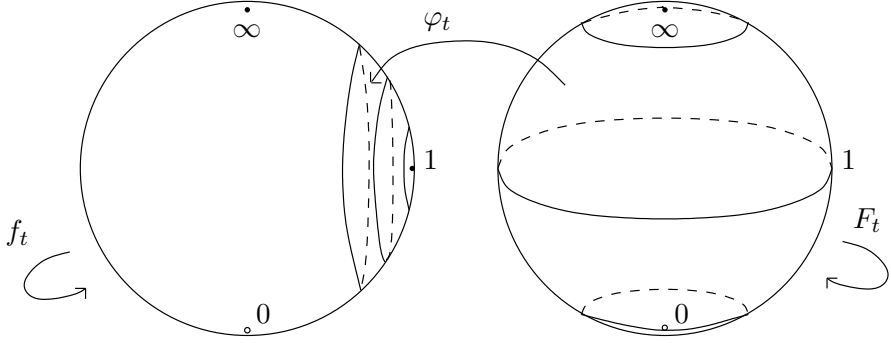}
\end{center}
\caption{Normal forms (\ref{Nf}) and (\ref{NF}) and transformation $\phi_t$ for small $t$.}
\label{figure3}
\end{figure}

Fix $a\in\hat{\C} -\{0,1,\infty\}$ and any continuous path $p:(0,1]\to M_2$ 
such that $p(t) \to [\Lambda_a]$ in $\Mbar_2$ as $t\to 0$.  For any representative of 
$p(t)$ in $\Rat_2$, the fixed points are distinct, so we can 
label the multipliers continuously so that $\alpha(t)\to a$, $\beta(t)\to 1/a$, and 
$\gamma(t)\to\infty$ as $t\to 0$.  Choose
a continuous path $\tilde{p}:(0,1]\to\Rat_2$ so that $[\tilde{p}] = p$
and $\tilde{p}(t)$ is normalized as in (\ref{Nf}). Then 
$\tilde{p}(t) \to \Lambda_a$ in $\Ratbar_2$ as $t\to 0$.  If we also label
the critical points of $\tilde{p}(t)$, then there is a unique M\"obius transformation
$\phi_t$ which transforms $\tilde{p}(t)$ into normalization (\ref{NF}).  
Then $\phi_t$ satisfies 
\begin{equation} \label{Nphi}
  \phi_t(z) = 1 + z\sqrt{\eps(t)} + o(\sqrt{\eps(t)}),
\end{equation}
for an appropriate choice of the square root of $\eps(t) := 1 - \alpha(t)\beta(t)$, 
locally uniformly for $z\in\C$ \cite[\S3]{Epstein:bounded}. 
See Figure \ref{figure3}.

\bigskip

\bigskip\noindent{\bf Iteration on $M_2$}.  
Fix $a\in \hat{\C} - \{0,1,\infty\}$, and 
let $\{f_t\in\Rat_2: t\in(0,1]\}$ be a continuous family
such that $f_t \to \Lambda_a$ in 
$\Ratbar_2$ as $t\to 0$.  Recall the definition of the indeterminacy
locus $I(2)$ in $\Ratbar_2$ given in \S\ref{ratd}.  By Theorem 
\ref{Ratindet}, $\Lambda_a \not\in I(2)$ implies 
that the iterates $f_t^n$ converge in 
$\Ratbar_{2^n}$ to 
\begin{equation}  \label{f_an}
  \Lambda_a^n = \left(a^nz \prod_{i=0}^{n-1} (z - w/a^i)^{2^{n-1-i}}:
                            w \prod_{i=0}^{n-1} (z - w/a^i)^{2^{n-1-i}}\right),
\end{equation}
by \cite[Lemma 7]{D:measures}, and therefore 
  $$f_t^n(z) \to a^n z$$
as $t\to 0$, 
locally uniformly on $\hat{\C} - \{1, 1/a, 1/a^2, \ldots, 1/a^{n-1}\}$
by Lemma \ref{sequence}.

Suppose further that $a = \zeta$  is 
a primitive $q$-th root of unity for some $q>1$.  
Then the family of $q$-th iterates $\{f_t^q\}$ converges to the 
identity function as $t\to 0$, 
locally uniformly on $\hat{\C} - \{1, \zeta, \ldots, \zeta^{q-1}\}$.  
However, Epstein showed that if we conjugate this family by a 
M\"obius transformation satisfying (\ref{Nphi}), then the limit
of the $q$-th iterate is a degree 2 map with a parabolic fixed
point:

\begin{prop} \label{Prop2}  \cite[Prop 2]{Epstein:bounded}
Let $\zeta$ be a primitive $q$-th root of unity for $q>1$, and 
suppose that $\{f_t: t\in (0,1]\}\subset\Rat_2$ is a continuous family 
normalized as in (\ref{Nf}) such that $\alpha(t)\to \zeta$ and $\beta(t)
\to 1/\zeta$ as $t\to 0$.  Suppose also that as $t\to 0$, 
  $$ \frac{\alpha(t)^q - 1}{\sqrt{\eps(t)}} \to \tau\in\Chat $$
for some choice of the square root of $\eps(t) = 1 - \alpha(t)\beta(t)$.  
Let $F_t = \phi_t^{-1} f_t \phi_t$ where $\phi_t$ satisfies (\ref{Nphi}) for 
this same choice of $\sqrt{\eps(t)}$.  Then  
$$  F_t^q(z) \to \left\{  \begin{array}{ll}
G_{\tau}(z)  = z + \tau + \frac{1}{z} &  \mbox{ for } \tau\in\C  \\
\infty & \mbox{ for } \tau=\infty
\end{array}  \right.  $$
as $t\to 0$, locally uniformly on $\C^*$.
\end{prop}

\bigskip\noindent
{\bf The $\tau^2$-value for embedded disks.}
In order to fully understand the iterate map, 
  $$\Phi_n: \Mbar_2\dashrightarrow \Mbar_{2^n},$$
defined by $[f]\mapsto [f^n]$ on $M_2$, 
we will need to analyse in more detail the behavior of the iterates
near the boundary point $\Lambda_\zeta\in\Ratbar_2$,
$\zeta^q = 1$.

Fix $q\geq 2$ and  $\zeta$ a primitive $q$-th root of unity.  
Let $\Delta: \D\hookrightarrow \Mbar_2$ be an embedded holomorphic
disk such that $\Delta(0) = [\Lambda_\zeta]$.    For $q\geq 3$ (so that 
$\zeta\not= 1/\zeta$), we can holomorphically parameterize two of 
the fixed point multipliers $\alpha(t) \to \zeta$ and $\beta(t) \to
1/\zeta$ as $t\to 0$ and set $\eps(t) = 1 - \alpha(t)\beta(t)$.  
We define the value $\tau^2\in\Chat$ for the disk $\Delta$ by
 $$ \tau^2(\Delta) = 
\lim_{t\to0} \frac{(\alpha(t)^q - 1)^2}{\eps(t)}$$
Set $\tau^2(\Delta) = \infty$ for any $\Delta\subset\del M_2$, 
since $\eps(t) \equiv 0$.  
If we interchange the labeling
of $\alpha$ and $\beta$, the $\tau^2$ value is unchanged: 
for any $q\geq 2$, the definition of $\eps = 1 - \alpha\beta$ implies that 
\begin{equation} \label{tau=}
\tau^2(\Delta) = \lim_{t\to0} \frac{(\alpha(t)^q - 1)^2}{\eps(t)} 
= \lim_{t\to0} \frac{(\beta(t)^q - 1)^2}{\eps(t)}.
\end{equation}
In the case of $q=2$, we are not necessarily able to holomorphically
label the fixed multipliers
$\alpha(t)$ and $\beta(t)$, since both tend to $-1 = \zeta$
as $t\to 0$.  Nevertheless, 
we can compute $\tau^2(\Delta)$ independent of 
any choice because of the equality in (\ref{tau=}).

\begin{theorem}  \label{tau2}
Let $\zeta$ be a primitive $q$-th root of unity for some $q\geq 2$, 
and let $\Delta$ and $\Delta'$ be two holomorphic disks in $\Mbar_2$
such that $\Delta(0) = \Delta'(0) = [\Lambda_\zeta]$.  Then for each $n\geq q$,
we have 
  $$\lim_{t\to 0} \Phi_n(\Delta(t)) = \lim_{t\to 0} \Phi_n(\Delta'(t))$$
in $\Mbar_{2^n}$ if and only if $\tau^2(\Delta) = \tau^2(\Delta')$.  
\end{theorem}

We will give the proof of Theorem \ref{tau2} in the next section.
For now, we complete the proof of Theorem \ref{indet}, which describes
the indeterminacy locus of the iterate map $\Phi_n$.

\bigskip\noindent {\bf Proof of Theorem \ref{indet}}.
Let $I(\Phi_n)$ denote the indeterminacy locus of the iterate map 
$\Phi_n: \Mbar_2 \dashrightarrow \Mbar_{2^n}$, and consider
the family $\Lambda_a\in\Ratbar_2$ for $a\in\Chat$ defined by 
(\ref{lambda}).  Since 
$\Mbar_2\iso\P^2$ is smooth, it suffices to show that 
$\Phi_n$ is discontinuous at $[f]\in\Mbar_2$ if and only 
if $[f] = [\Lambda_\zeta]$ for a primitive $q$-th root of 
unity $\zeta$, with $1 < q \leq n$.  

Suppose $a\in\hat{\C}$ is not a primitive $q$-th root 
of unity.  By Lemma \ref{stableLambda}, $\Lambda_a^n\in\Ratbar_{2^n}$
is stable for all $n\geq 1$, so that by Lemma \ref{stableindet}, 
$[\Lambda_a] \not\in I(\Phi_n)$ for all $n\geq 1$.  

Fix $q\geq 2$ and let $\zeta$ be a primitive $q$-th 
root of unity.   By Lemma 
\ref{stableLambda}, the iterate $\Lambda_\zeta^n\in\Ratbar_{2^n}$ is stable 
if and only if 
$1 \leq n< q$, and so by Lemma \ref{stableindet}, 
$[\Lambda_\zeta]\not\in I(\Phi_n)$ for all $n<q$.  

Now suppose $n\geq q$.  By Theorem \ref{tau2}, it suffices
to show there exist holomorphic disks $\Delta$ and $\Delta'$ in $\Mbar_2$
such that $\Delta(0) =\Delta'(0) = [\Lambda_\zeta]$ 
and $\tau^2(\Delta)\not= \tau^2(\Delta')$.  
Define $f_t\in\Rat_2$ by normal form (\ref{Nf}) with 
$\alpha(t) = \zeta + t$ and $\beta(t) = 1/\zeta + 2 t$.  Then it 
is easy to compute that $\tau^2(\Delta) = 0$ for $\Delta(t) = [f_t]$.
On the other hand, 
$\tau^2(\Delta') = \infty$ for any disk $\Delta' \subset \del M_2$.
\qed

\bigskip
\section{The iterate map in degree 2}
\label{degree2proofs}

Fix $q\geq 2$ and let $\zeta$ be a primitive $q$-th root of unity.  
Let $\Lambda_\zeta\in\Ratbar_2$ be defined by equation (\ref{lambda}).
In this section, we prove Theorem \ref{tau2}, that the limiting values 
of the iterate map $\Phi_n: \Mbar_2\dashrightarrow\Mbar_{2^n}$ 
on a holomorphic disk passing through $[\Lambda_\zeta]$ 
depends only on the $\tau^2$-value of the disk.  
We treat the 
cases of $\tau^2\in\C$ and $\tau^2=\infty$ separately (Propositions  
\ref{Fiterates} and \ref{paraboliclimit}).

\begin{prop} \label{Fiterates}  Fix $q>1$ and $\zeta$ a primitive 
$q$-th root of unity.
Suppose that $\{f_t: t\in (0,1]\}\subset\Rat_2$ is a continuous family, 
normalized as in (\ref{Nf}), such that
$\alpha(t) \to \zeta$, $\beta(t) \to 1/\zeta$,
and $(\alpha(t)^q - 1)/\sqrt{\eps(t)} \to \tau\in\C$,
for some choice of $\sqrt{\eps(t)}$ as $t\to 0$.  Let $\phi_t \in\Aut\Chat$
satisfy (\ref{Nphi}) for this choice of $\sqrt{\eps(t)}$.  
Then the the $n$-th iterate of $F_t = \phi_t^{-1}f_t\phi_t$ converges
in $\Ratbar_{2^n}$ as $t\to 0$ to the following:  
  $$F_{q,\tau,n} =  \left\{ \begin{array}{ll}
     (z^{2^{n-1}}w^{2^{n-1}}:0) &   \mbox{ for } 1\leq n  < q, \\
     (z^{2^{q-1}-1}w^{2^{q-1}-1}(z^2 + \tau zw + w^2):
                      z^{2^{q-1}}w^{2^{q-1}}) &  \mbox{ for } n=q, \\
    F_{q,\tau,(n\mod q)} \circ (F_{q,\tau,q})^{\lfloor n/q \rfloor}
                     &  \mbox{ for } n > q,   
                     \end{array}   \right.  $$
where $F_{q,\tau,0}$ is the identity map.  
\end{prop}

\noindent
Recall that by Proposition \ref{Prop2}, the $q$-th iterate 
of $F_t$ converges to $G_\tau(z) = z + \tau + 1/z$, locally
uniformly on $\C^*$.  By Proposition \ref{Fiterates} 
(together with Lemma \ref{sequence}), 
$F_t^n(z) \to \infty$ as $t\to 0$
locally uniformly on the complement of a finite set 
in $\hat{\C}$ for all $n$ which are not multiples of $q$,
and for every $n=mq$, $F_t^{mq}(z) \to G_{\tau}^m(z)$ as $t\to 0$
locally uniformly on the complement of a finite 
set in $\hat{\C}$.

\begin{lemma} \label{17}  \cite[\S4, (17)]{Epstein:bounded}
Fix $a\in\Chat - \{0,1,\infty\}$ and 
suppose that $\{f_t: t\in (0,1]\}\subset\Rat_2$ is a continuous family
normalized as in (\ref{Nf}) such that $f_t \to \Lambda_a$ in $\Ratbar_2$ 
as $t\to 0$.  Let $\eps(t) = 1 - \alpha(t)\beta(t)$ and let 
$z: (0,1]\to \Chat$ be a continuous path.  Then 
$$  \frac{f_t(z(t))}{z(t)} = \left\{  \begin{array}{ll}
\alpha(t) + o(1) & \mbox{ if } \eps(t) = o(z(t)-1)  \\
\alpha(t) + o(\sqrt{\eps(t)}) & \mbox{ if } \sqrt{\eps(t)} = o(z(t)-1)
\end{array} \right.  $$
as $t\to 0$.
\end{lemma}

\bigskip\noindent {\bf Proof of Proposition \ref{Fiterates}}.
For each fixed $z\in\C^*$, 
we have $\eps(t) = o(\phi_t(z) - 1)$ and therefore, 
by Lemma \ref{17}, 
  $$f_t(\phi_t(z))/\phi_t(z) = \alpha(t) + o(1).$$
In particular, $f_t(\phi_t(z)) \to \zeta,$
locally uniformly on $\C^*$.  Since $\zeta\not= 1$, 
we obtain 
  $$F_t(z) = \phi_t^{-1} f_t \phi_t(z) \to \infty,$$
locally uniformly in $\C^*$.  By induction we find that 
$\eps(t) = o(f_t^{n-1}(\phi_t(z)) - 1)$ for each $1\leq n \leq q$, so that 
\begin{equation} \label{limfn}
  f_t^n(\phi_t(z)) = \alpha(t) f_t^{n-1}(\phi_t(z)) + o(1) \to \zeta^n 
\end{equation} 
as $t\to 0$.  Consequently, $F_t^n(z) \to \infty$ locally 
uniformly in $\C^*$ for each $n< q$.  

It follows that for each $n< q$, the limit of 
$F_t^n$ exists in $\Ratbar_{2^n}$ as $t\to 0$:  by Lemma \ref{sequence}(ii)
it must be of the form $F_{q,\tau,n} = (z^kw^l: 0)$ for non-negative
integers $k$ and $l$ such that $k+l = 2^n$ because the convergence
of $F_t^n(z) \to \infty$ is uniform away from 0 and $\infty$.
To determine $k$ and $l$, 
it suffices (again by Lemma \ref{sequence}) to count the preimages of 0
by $F_t^n$ near both 0 and $\infty$.

Fix $n< q$.  The iterate $f_t^n \to \Lambda_\zeta^n$ in $\Ratbar_{2^n}$
as $t\to 0$, where $\Lambda_\zeta^n$ is given 
by equation (\ref{f_an}), and Lemma \ref{sequence} implies that
$f_t^n(z) \to \zeta^n z$ as $t\to 0$, 
locally uniformly on $\Chat - \{1, 1/\zeta, \dots, 1/\zeta^{n-1}\}$.
Therefore, for all sufficiently small $t$, there is exactly one preimage by 
$f_t^n$ of $z=1$ very close to $z = 1/\zeta^n$.   
Fix small disks $\overline{D}_1 \subset D_2$ around $z=1$.
Counting the depths of the holes outside the disk $D_1$, Lemma
\ref{sequence} implies that there are exactly $2^{n-1}$ preimages
of $z=1$ by $f_t^n$ in $\Chat - D_1$ for sufficiently small $t$.  
Thus, for $F_t^n = \phi_t^{-1}f_t^n \phi_t$, there are exactly $2^{n-1}$ preimages of
0 in $\phi_t^{-1}(\Chat-D_1)$ for sufficiently small $t$.  Therefore
the depth of $z=\infty$ for $F_{q,\tau,n}$ is at 
least $2^{n-1}$.  

On the other hand, 
let $D_0$ be any disk around $z=0$.  As in the argument 
to show  (\ref{limfn}), 
$f_t^n(\overline{D}_1 - \phi_t(D_0)) \subset \zeta^n D_2$
for all sufficiently small $t$ (and $n<q$).    Therefore, 
$F_t^n \to \infty$ uniformly on $\phi_t^{-1} (D_1) - D_0$,
so that $F_t^n$ has at most $2^{n-1}$ preimages of 0 
outside $D_0$ for small $t$.  Therefore, 
the depth of $F_{q,\tau,n}$ at $\infty$ is exactly $2^{n-1}$, 
and we can conclude that 
  $$F_{q,\tau,n} = (z^{2^{n-1}}w^{2^{n-1}}: 0).$$
  
Now suppose $n=q$.  Since $F_t^q(z) \to G_{\tau}(z) = 
z + \tau + 1/z$ locally 
uniformly in $\C^*$ by Proposition \ref{Prop2}, 
the limit of $F_t^q$ must exist in $\Ratbar_{2^q}$ and 
be of the form $F_{q,\tau,q} = (z^kw^l(z^2+\tau zw + w^2):z^{k+1}w^{l+1})$
for integers $k$ and $l$ with $k+l = 2^q-2$.  To compute 
$k$ and $l$, 
we will count preimages of $z=\infty$ near $\infty$.

Without loss of generality, we may assume that $\phi_t$
fixes $\infty$ for all $t$. 
As before, let $D_0$ be a disk centered at $z=0$ and 
let $D_1$ be a small disk around $z=1$.  Because of 
the depths of the holes of $\Lambda_\zeta^q$, we find
that there are exactly $2^{q-1}$ preimages of $\infty$
by $f_t^q$ in $\Chat- D_1$ for all sufficiently small $t$ 
($2^{q-1}-1$ of them 
accumulate on the $q$-th roots of unity and one preimage
is at $\infty$). 
Also, from (\ref{limfn}) when $n=q$, there are no preimages
of $\infty$ by $f_t^q$ in $D_1-\phi_t(D_0)$.   Therefore, 
$F_t^q$ has exactly $2^{q-1}$ preimages of $\infty$ in 
$\Chat - D_0$ for all sufficiently small $t$.  Consequently, 
$F_{q,\tau,q}$ has a hole of depth exactly $2^{q-1}-1$
at $\infty$, and therefore, 
  $$F_{q,\tau,q} = (z^{2^{q-1}-1}w^{2^{q-1}-1}
         (z^2 + \tau zw + w^2):  z^{2^{q-1}}w^{2^{q-1}}).$$

Finally we need to compute the limits of $F_t^n$ in $\Ratbar_{2^n}$ for 
every $n > q$.  Write $n = k+mq$ for integers $0\leq k<q$ and $m>0$.  
For $k=0$ the 
desired form follows immediately from Theorem \ref{Ratindet} since $F_{q,\tau,q}
\not\in I(2^q)$.  For $k>0$, the result follows from 
Lemma \ref{composition}, since 
$F_t^n = F_t^k\circ F_t^{mq}$.  
\qed

\medskip
\begin{lemma} \label{Fmeasures}
For $q\geq 2$ and $\tau\in\C$, let $\{F_t\}$ be a family of rational maps as in
Proposition \ref{Fiterates}.  
Then the measures $\mu_{F_t}$ converge weakly
to $\mu_{F_{q,\tau,q}}$ as $t\to 0$.  Furthermore, 
$\mu_{F_{q,\tau,q}}(\{z\}) < 1/2$
for all $z\in\Chat$. 
\end{lemma}

\proof
The convergence of the measures follows immediately
from Theorem \ref{limitmu} because $F_{q,\tau,q} \not\in I(2^q)$.
Let $\mu_\tau = \mu_{F_{q,\tau,q}}$ as defined in Section \ref{ratd}.  
By Lemma \ref{Lemma5}, we can compute the values of $\mu_\tau$, 
 $$\mu_\tau(\{\infty\}) = \mu_\tau(\{0\}) 
= \frac{2^{q-1} -1}{2^q} \sum_{l=0}^\infty \frac{1}{2^{ql}}
= \frac{2^{q-1} - 1}{2^q-1} < \frac12.$$
Then, for any point $p\in\C^*$, we have 
  $$\mu_\tau(\{p\}) \leq 1 - \mu_\tau(\{0\}) - \mu_\tau(\{\infty\}) 
= \frac{1}{2^q - 1} < \frac12.$$
\qed

\begin{prop} \label{paraboliclimit}
Fix $q>1$ and $\zeta$ a primitive $q$-th root of unity.  
Suppose that $\{f_t: t\in(0,1]\}\subset\Rat_2$ is a continuous family of rational
maps normalized as in (\ref{Nf}), such that  
$\alpha(t) \to \zeta$, $\beta(t) \to 1/\zeta$, and 
$(\alpha(t)^q - 1)/\sqrt{\eps(t)} \to \infty$ as $t\to 0$.  Conjugating
by $A_t(z) = 1 + z(\alpha(t)^q-1) \in \Aut\Chat$, the 
iterates of $P_t = A_t^{-1} f_t A_t$ converge
in $\Ratbar_{2^n}$ as $t\to 0$ to the following:  
  $$P_{q,n} =  \left\{ \begin{array}{ll}
     (z^{2^{n-1}}w^{2^{n-1}}:0) &   \mbox{ for } 1\leq n  < q, \\
     (z^{2^{q-1}}w^{2^{q-1}-1}(z+w):
                      z^{2^{q-1}}w^{2^{q-1}}) &  \mbox{ for } n=q, \\
    P_{q,(n\mod q)} \circ (P_{q,q})^{\lfloor n/q \rfloor}, 
                     &  \mbox{ for } n > q,   
                     \end{array}   \right.  $$
where $P_{q,0}$ is the identity map.  
\end{prop}

\noindent
In particular, Lemma \ref{sequence} implies that
$P_t^n(z) = A_t^{-1}f_t^nA_t(z) \to \infty$ as $t\to 0$, 
locally uniformly on the complement of a finite set 
in $\hat{\C}$,  for all $n$ which are not multiples of $q$,
and $P_t^{mq}(z) \to z+m$ as $t\to 0$,
locally uniformly on $\hat{\C}-\{0,-1,-2, \ldots, -m+1, \infty\}$ 
for every $m\geq 1$.  

\proof
The proof is similar to the proof of Proposition \ref{Fiterates}.
For fixed $z\in\C^*$,
we have $A_t(z) - 1 = z(\alpha(t)^q - 1)$, and so 
$\sqrt{\eps(t)} = o(A_t(z)-1)$.  By Lemma \ref{17},
  $$\frac{f_t(A_t(z))}{A_t(z)} = \alpha(t) + o(\sqrt{\eps(t)}).$$
In particular, $f_t(A_t(z)) \to \zeta$, locally uniformly in $\C^*$.
Thus $\sqrt{\eps(t)} = o(f_t(A_t(z)) - 1)$ also and so 
  $$\frac{f_t^2(A_t(z))}{f_t(A_t(z))} = \alpha(t) + o(\sqrt{\eps(t)}).$$
By induction, we have 
\begin{equation}  \label{fqA}
 \frac{f_t^q(A_t(z))}{A_t(z)} = \prod_1^q \frac{f_t^n(A_t(z))}{f_t^{n-1}(A_t(z))}
       = \alpha(t)^q + o(\sqrt{\eps(t)}),
\end{equation}
and we conclude that 
  $$P_t^q(z) = A_t^{-1} f_t^q A_t(z) = \frac{\alpha(t)^q-1}{\alpha(t)^q-1} 
     + \frac{\alpha(t)^q(\alpha(t)^q -1) z}{\alpha(t)^q - 1} + 
      \frac{o(\sqrt{\eps(t)})}{\alpha(t)^q - 1}  \to z+ 1$$
as $t\to 0$, and that $P_t^n(z) \to \infty$ for each $n<q$ as 
$t\to 0$, locally uniformly on $\C^*$.  

It remains to determine the limit of the iterates of $P_t$ in
$\Ratbar_{2^n}$.  Fix $n<q$.  The proof that $P_t^n \to
P_{q,n}$ in $\Ratbar_{2^n}$ is identical to the proof of 
Proposition \ref{Fiterates} in the case $n<q$, and we omit it.  

Let $n=q$.  Since $P_t^q(z) \to z+1$ locally uniformly on 
$\C^*$, Lemma \ref{sequence} implies that the limit of $P_t^q$ 
must exist in $\Ratbar_{2^q}$ and be of the form $P_{q,q} = 
(z^kw^l(z+w):z^kw^{l+1})$
for integers $k$ and $l$ such that $k+l = 2^q-1$.  
As in the proof of Proposition \ref{Fiterates}, the convergence
of $f_t^q \to \Lambda_\zeta^q\in \Ratbar_{2^q}$ and the
estimate on $f_t^q(A_t(z))$ in (\ref{fqA}) imply that 
$P_t^q$ will have exactly $2^{q-1}$ preimages of $\infty$
near $\infty$ for all sufficiently small $t$.  Therefore, $l+1 = 
2^{q-1}$ and $P_{q,q}$ will have a hole of depth $2^{q-1}-1$
at $\infty$.  

For each $n>q$, the formula for $P_{q,n}$ follows from 
Lemma \ref{composition}.
\qed

\bigskip\noindent{\bf Stability of the iterates.}
Recall the criteria for GIT stability of points in $\Ratbar_{2^n}$ given
in Section \ref{GITstability}.  For each $n\geq 2$, define $F_{q,\tau,n}$
and $P_{q,n}$ in $\Ratbar_{2^n}$ as in Propositions 
\ref{Fiterates} and \ref{paraboliclimit}.  The following lemma shows
that $F_{q,\tau,n}$ and $P_{q,n}$ define points in $\Mbar_{2^n}$ 
if and only if $n\geq q$.

\begin{lemma} \label{stablelimit}
Let $\zeta$ be a primitive $q$-th root of unity for $q\geq 2$ 
and fix $\tau \in\C$.  Then  each of $F_{q,\tau, n}$ 
and $P_{q,n} \in
\Ratbar_{2^n}$ is stable if and only if $n\geq q$.
\end{lemma} 

\proof
For $n<q$, $F_{q,\tau,n} = P_{q,n} = (z^{2^{n-1}}w^{2^{n-1}}:0)$
has a hole of depth $2^{n-1} = 2^n/2$ at $\infty$, and the associated
rational map of lower degree is the constant $\infty$.  Thus, $F_{q,\tau,n}
= P_{q,n}$ is not stable.  

Fix $\tau\in\C$.  Note that $F_{q,\tau,q}\not\in I(2^q)$, so it has well-defined
forward iterates.  By Lemma \ref{Fmeasures}, the measure 
$\mu_\tau = \mu_{F_{q,\tau,q}}$ has no atoms of mass $\geq 1/2$, 
and therefore $F_{q,\tau,q}$ and all its forward iterates 
$F_{q,\tau,qm}$ are stable by Proposition \ref{evendeg}.   

Fix integers $0<k<q$ and $m\geq 1$.  
Lemma \ref{Cor8} allows us to compute the 
depths of 0 and $\infty$ as holes of $F_{q,\tau,qm}$:  
$$d_0(F_{q,\tau,qm}) = d_\infty(F_{q,\tau,qm})
= 2^{qm} \left(\frac{2^{q-1}-1}{2^q} \sum_{l=0}^{m-1} \frac{1}{2^{ql}} \right) 
= \frac{2^{q-1}-1}{2^q-1} \left( 2^{qm} - 1 \right).$$ 
With Lemma \ref{composition} we 
can compute the depths of the holes of $F_{q,\tau,k+mq}$:
$$ d_0(F_{q,\tau,k+mq}) = d_\infty(F_{q,\tau,k+mq}) 
= 2^k \frac{2^{q-1}-1}{2^q-1} (2^{qm} - 1) + 2^{k-1} 
\qquad\qquad\qquad $$  $$\qquad\qquad\qquad
= 2^{k+qm} \left( \frac{2^{q-1}-1}{2^q-1} + \frac{1}{2^{qm}}
\left( \frac12 - \frac{2^{q-1}-1}{2^q-1}\right) \right)
< \frac{2^{k+qm}}{2}, $$
and therefore for any $z\not= 0,\infty$, 
$$ d_z(F_{q,\tau,k+mq}) \leq 2^{k+mq} - 2 d_0(F_{q,\tau,k+mq}) 
\leq 2^{k+mq}/ (2^q - 1).  $$
This implies that $F_{q,\tau,k+mq}$ is stable.  

Now consider 
$P_{q,q} = (z^{2^{q-1}}w^{2^{q-1}-1} (z+w): z^{2^{q-1}}w^{2^{q-1}})$.
By the definition of the measure $\mu_{P_{q,q}}$, 
  $$\mu_{P_{q,q}} = \frac{1}{2^q} \sum_{k=0}^\infty \frac{2^{q-1}}{2^{qk}}
\delta_{-k}  +  \frac{1}{2^q} \sum \frac{2^{q-1}-1}{2^{qk}} \delta_\infty 
= \frac12 \sum_{k=0}^\infty \frac{1}{2^{qk}} \delta_{-k} 
+ \frac{2^{q-1}-1}{2^q-1} \delta_\infty.$$
By Proposition \ref{evendeg}, $P_{mq}$ is stable for all $m\geq 1$. 

Fix integers $0<k<q$ and $m\geq 1$.  
Lemma \ref{Cor8} allows us to compute the depths of the holes 
of $P_{q,qm}$:
$d_0(P_{q,qm}) = 2^{qm-1}$, $d_z(P_{q,qm}) \leq 2^{qm-2}$
for all $z\not= 0,\infty$, and 
 $$d_\infty(P_{q,qm}) = \frac{2^{qm}}{2^q} \sum_{l=0}^{m-1} 
\frac{2^{q-1}-1}{2^{ql}} = \frac{2^{q-1}-1}{2^q-1} (2^{qm} - 1).$$
Therefore, by Lemma \ref{composition}, 
 $$d_\infty(P_{q,k+qm}) = 2^k \frac{2^{q-1}-1}{2^q-1} (2^{qm} - 1) + 2^{k-1}
< 2^{k+qm}/2,$$
$d_0(P_{q,k+qm}) = 2^{k+qm}/2$, and $d_z(P_{q,k+qm}) \leq 
2^{q+km}(1/4 + 1/2^{qm+1}) < 2^{k+qm}/2$ for all $z\not= 0,\infty$.
Observing that $\phi_{P_{q,k+qm}}(0)\not=0$, we see that $P_{q,k+qm}$
is stable. 
\qed

\bigskip\noindent{\bf A non-constant map from $\Chat$ to the boundary
of $M_{2^n}$.}
Recall the definitions of $F_{q,\tau,n}$ and $P_{q,n}$ from Propositions
\ref{Fiterates} and \ref{paraboliclimit}.  
The next two lemmas show that there are regular maps of 
degree 2 from $\P^1$ to the boundary of $M_{2^n}$ in $\Mbar_{2^n}$
parameterized by the families $F_{q,\tau,n}$ for each $q\leq n$. 

\begin{lemma} \label{nonconstant}
Fix integers $n\geq q\geq 2$ and $\sigma, \tau\in\C$.  The following are 
equivalent:  
\begin{itemize}
\item[(i)]  $[F_{q,\tau,n}] = [F_{q,\sigma,n}]$ in $\Mbar_{2^n}$,
\item[(ii)]  $\mu_{F_{q,\tau,q}} = A_* \mu_{F_{q,\sigma,q}}$ for some 
        $A\in \Aut\Chat$, and 
\item[(iii)] $\sigma = \pm \tau$. 
\end{itemize}
Furthermore,  $[P_{q,n}] \not= [F_{q,\tau,n}]$ for all $\tau\in\C$.  
\end{lemma}

\proof
Write $n = k + qm$ for integers $0\leq k < q$ and $m>0$, and let $\mu_\tau 
= \mu_{F_{q,\tau,q}}$.    
Define $A\in\Aut\Chat$ by 
$A(z:w) = (-z:w)$, and note that 
$A F_{q,\tau,q} A^{-1} = F_{q,-\tau,q}$.
Therefore, 
 $$A F_{q,\tau,qm} A^{-1} = A (F_{q,\tau,q})^m A^{-1} = F_{q,-\tau,qm},$$
so that $[F_{q,\tau,qm}] = [F_{q,-\tau,qm}]$ and $\mu_\tau = A_*\mu_{-\tau}$.
For $k\geq 1$, note that $A F_{q,\tau,k} A^{-1} = (\pm z^{2^{k-1}}w^{2^{k-1}}: 0) = 
F_{q,\tau,k}$ and is independent of $\tau$, so that 
 $$A F_{q,\tau,k+qm} A^{-1} = A F_{q,\tau,k} F_{q,\tau,mq} A^{-1}
  = F_{q,\tau,k} F_{q,-\tau,qm} = F_{q,-\tau,k+qm},$$
and therefore, $[F_{q,\tau,n}] = [F_{q,-\tau,n}]$.  
This proves that (iii) implies both (i) and (ii).  

To see that (i) implies (iii), note that for each $\tau\in\C$ and $n\geq q$, 
the points 0 and $\infty$ in $\Chat$ are distinguished
by their depths as holes of $F_{q,\tau, n}$:  as computed in the previous lemma, 
the depth at 0 is the same as the depth at $\infty$ and greater 
than at any other point.  
The two preimages of 0 by the 
degree two map $z + \tau + 1/z$ are also distinguished:  
if $n=q$, then they are sent to 0 by the dynamics, and if $n>q$, 
they are distinguished by their depths which can be computed 
with Lemma \ref{composition}.  
Therefore, if $F_{q,\tau,n}$ and $F_{q,\sigma,n}$ are equivalent,
the cross-ratio of these four 
points must coincide.  

The preimages of 0 by $z+ \tau + 1/z$ lie at the points 
  $$p_{\pm} := \frac12 (-\tau \pm \sqrt{ \tau^2 - 4}).$$
We compute the cross ratio $\chi$, normalized so that
  $$\chi(0,\infty,1,z) = z.$$
Then 
  $$\chi(\tau) := \chi(0,\infty, p_+, p_-) = p_-/p_+ 
= \frac{1}{4} \left(\tau + \sqrt{\tau^2-4}\right)^2.$$
Notice that reversing the labeling of $p_+$ and $p_-$ 
or of 0 and $\infty$ gives $1/\chi(\tau) = \chi(-\tau)$, and
so the natural invariant to consider is 
  $$\chi(\tau) + 1/\chi(\tau),$$
and this is what we will compute.  We find, 
  $$\chi + 1/\chi = \tau^2 - 2.$$
Therefore, if $F_{q,\tau,n}$ and $F_{q,\sigma,n}$
 are equivalent,
we must have $\tau^2 -2 = \sigma^2 -2$, and therefore, 
$\sigma = \pm \tau$.  

The proof that (ii) implies (iii) is similar:  the (unordered) pairs of points 
$\{0, \infty\}$ and $\{p_+,p_-\}$ are distinguished by their masses, as can be computed
with Lemma \ref{Lemma5}.  If (ii) holds, the cross ratio of these
four points must coincide for $\tau$ and for $\sigma$.  As seen above,
this implies that $\sigma = \pm \tau$.

Finally, in the proof Lemma \ref{stablelimit} we showed that 
$d_0(P_{q,n}) = 2^{n-1}$ which is strictly greater than the
depth of any hole for $F_{q,\tau,n}$.  Therefore, 
$[F_{q,\tau,n}] \not= [P_{q,n}]$.
\qed

\bigskip
\begin{lemma} \label{continuous}
For each $n\geq q\geq 2$, $[F_{q,\tau,n}] \to [P_{q,n}]$ 
in $\Mbar_{2^n}$ as $\tau\to \infty$.
\end{lemma}

\proof
Write $n = k + mq$ for integers $0\leq k < q$ and $m\geq 1$.  
Recall the definitions,  
 $$F_{q,\tau,q} = (z^{2^{q-1}-1}w^{2^{q-1}-1}(z^2 + \tau zw + w^2): 
 z^{2^{q-1}}w^{2^{q-1}}) \in\Ratbar_{2^q},$$
and
  $$P_{q,q} = (z^{2^{q-1}}w^{2^{q-1}-1}(z+w): z^{2^{q-1}}w^{2^{q-1}})
  \in\Ratbar_{2^q}.$$
For each $\tau\in\C$, define $A_\tau\in\Aut\Chat$ by 
  $$A_\tau(z:w) = (\tau^{-1}(z+w): w).$$
Then 
  $$A_\tau F_{q,\tau,q} A_\tau^{-1} 
  = (\tau^{2^{q-1}} z^{2^{q-1}-1}w^{2^{q-1}-1}(z^2 + zw) + O(\tau^{2^{q-1}-1}):
  \tau^{2^{q-1}} z^{2^{q-1}}w^{2^{q-1}} ),$$
and therefore $A_\tau F_{q,\tau,q} A_\tau^{-1} \to P_{q,q}$ in $\Ratbar_{2^q}$
as $\tau\to \infty$.  By the regularity of the iterate maps near $P_{q,q}$, we have
also that $A_\tau F_{q,\tau,q}^m A_\tau^{-1} \to P_{q,q}^m = P_{q,qm}$ 
in $\Ratbar_{2^{qm}}$ as $\tau\to \infty$.  Note also that for $1\leq k < q$, 
  $$A_\tau F_{q,\tau,k} A_\tau^{-1} = (\tau^{-1} (\tau z - w)^{2^{k-1}} w^{2^{k-1}}: 0) 
  \to P_{q,k} = (z^{2^{k-1}}w^{2^{k-1}}:0)$$
in $\Ratbar_{2^k}$ as $\tau\to \infty$.  Therefore, by the continuity of the 
composition map (Lemma \ref{composition}), 
  $$A_\tau F_{q,\tau,n} A_\tau^{-1} = A_\tau F_{q,\tau,k} F_{q,\tau,qm} A_\tau^{-1}
  \to P_{q,k} P_{q,qm} = P_{q,n},$$
in $\Ratbar_{2^n}$ as $\tau\to \infty$.  Since each of these elements 
in $\Ratbar_{2^n}$ is stable by Lemma 
\ref{stablelimit}, we can conclude that $[F_{q,\tau,n}] \to [P_{q,n}]$ in 
$\Mbar_{2^n}$ as $\tau\to \infty$.
\qed

\bigskip\noindent{\bf The iterate map $\Phi_n$.}
We are now ready to prove Theorem \ref{tau2}.  However, with all 
the notation of this section in place, we can state a more specific
result about the limiting values of the iterate map $\Phi_n: 
M_2 \to M_{2^n}$.  Recall the definition of the $\tau^2$-value
of a holomorphic disk in $\Mbar_2$, given just before the statement
of Theorem \ref{tau2}.  For integers $n\geq q\geq 2$ and 
$\tau\in\C$, let $F_{q,\tau,n}\in\Ratbar_{2^n}$ be 
defined as in Proposition \ref{Fiterates} and 
let $P_{q,n}\in\Ratbar_{2^n}$ be defined 
as in Proposition \ref{paraboliclimit}.  

\begin{prop}  \label{limtau}
Fix integers $n\geq q\geq 2$ and let $\zeta$ be a primitive 
$q$-th root of unity.  
Let $\Delta: \D\hookrightarrow \Mbar_2$ be a holomorphic 
disk such that $\Delta(0) = [\Lambda_\zeta]$.  
If $\tau^2(\Delta) = \infty$, then $\Phi_n(\Delta(t)) \to [P_{q,n}]$
in $\Mbar_{2^n}$ as $t\to 0$.  
If $\tau$ is a square root of $\tau^2(\Delta)\in\C$,
then $\Phi_n(\Delta(t)) \to [F_{q,\tau,n}]$  as $t\to 0$.  
\end{prop}

\proof
Suppose that $\Delta\subset
\del M_2$.  Then by definition, $\tau^2(\Delta) = \infty$.  
For any real path $p: [0,1] \to \D$ such that $p(0)=0$, choose a 
continuous lift of $\Delta(p(t))$ to $\Lambda_{a(t)}\in\Ratbar_2$
so that $\Delta(p(t)) = [\Lambda_{a(t)}]$ and $a(t)\to \zeta$ as 
$t\to 0$.  
Define $A_t\in\Aut\Chat$ 
by $A_t(z,w) = ((a(t)^q-1)z + w: w)$.  Using formula 
(\ref{f_an}) for the iterates of $\Lambda_{a(t)}$, 
it can be computed directly that $P_t^n := A_t^{-1} \Lambda_{a(t)}^n A_t 
\in\Ratbar_{2^n}$ is given by 
$$\begin{array}{l}
 \left(z^{2^{n-1}} (a(t)^n (a(t)^q - 1)z + (a(t)^n-1)w) 
  \left(\prod_{i=1}^{n-1} ( (a(t)^q-1)z + w(1 - 1/a(t)^i))^{2^{n-1-i}} \right): \right.  \\
     \left. \qquad  (a(t)^q-1) z^{2^{n-1}}w 
    \left(\prod_{i=1}^{n-1} ( (a(t)^q-1)z + w(1 - 1/a(t)^i))^{2^{n-1-i}} \right) \right).
 \end{array}$$
For each $n\leq q$, we see immediately that $P_t^n$ converges in 
$\Ratbar_{2^n}$ to $P_{q,n}$ as $t\to 0$,
where $P_{q,n}$ is defined in Proposition \ref{paraboliclimit}.
Lemma \ref{composition} then implies that $P_t^n \to P_{q,n}$
as $t\to 0$ for all $n\geq q$.   By Lemma \ref{stablelimit}, the 
point $P_{q,n}\in\Ratbar_{2^n}$ is stable (in the sense of GIT)
and so defines a conjugacy class in $\Mbar_{2^n}$.  
Consequently, 
 $$\lim_{t\to 0} \Phi_n(\Delta(t)) = [P_{q,n}].$$
If $\Delta\not\subset \del M_2$, but $\tau^2(\Delta) = \infty$,
then by Proposition \ref{paraboliclimit}, 
 $$\lim_{t\to 0} \Phi_n(\Delta(t)) = [P_{q,n}].$$
For $\tau^2(\Delta) \in\C$, let $\tau\in\C$ be a square root of  
$\tau^2(\Delta)$.  Then by Proposition \ref{Fiterates},
  $$\lim_{t\to 0} \Phi_n(\Delta(t)) = [F_{q,\tau,n}].$$   
\qed

\bigskip\noindent {\bf Proof of Theorem \ref{tau2}.}
By Lemma \ref{nonconstant}, the conjugacy classes 
$[P_{q,n}]$ and $[F_{q,\tau,n}]$ are all distinct in 
$\Mbar_{2^n}$ when $n\geq q$.  Therefore, 
the theorem is an immediate corollary of Proposition 
\ref{limtau}.
\qed

\bigskip

\section{The blow-ups of $\Mbar_2$}
\label{blowup}

Let $\Phi_n$ denote the $n$-th iterate map $M_2\to M_{2^n}$, and 
let $\Gamma_n$ denote the closure of the image of $M_2$ in the product
in $\Mbar_2\times\Mbar_4 \times \cdots \Mbar_{2^n}$ via the embedding
$(\Id, \Phi_2, \ldots, \Phi_n)$.  
Let $\pi_n: \Gamma_n\to\Gamma_{n-1}$ be
the projection to the first $n-1$ factors.  In this section we study the structure
of the pair $(\Gamma_n, \pi_n)$ for each $n\geq 2$, and we  
give the proof of Theorem \ref{finite}.  

Note that the composition $\pi_2\pi_3\cdots\pi_n:\Gamma_n\to \Mbar_2$ is an isomorphism away 
from the finite indeterminacy set $I(\Phi_n)\subset\Mbar_2$, 
described explicitly in Theorem \ref{indet}.  
The projection of $\Gamma_n$ to the $n$-th
factor $\Mbar_{2^n}$ is a regular extension of $\Phi_n$.  

\bigskip\noindent{\bf The model for $\Gamma_n$.} 
First, let $p_2: B_2 \to \Mbar_2\iso \P^2$ denote the standard blow-up of $\P^2$ 
at the unique point in $I(\Phi_2)$.  That is, $B_2$ is the closure of the graph of 
$\P^2 \dashrightarrow \P^1$ given by $(x:y:z) \mapsto (x:y)$, with the 
coordinates chosen so that $I(\Phi_2) = \{(0:0:1)\}$, and $p_2$ is the 
projection to the first factor.  Inductively define $p_n:
B_n \to B_{n-1}$
to be a blow-up of $B_{n-1}$ at each point of $I(\Phi_n) - I(\Phi_{n-1})$, 
which in local coordinates is given by the blow-up of $\C^2$ along the ideal
$(x^2,y)$ where the axis $\{y=0\}$ represents the boundary of $M_2$.  That is, 
$B_n$ is locally isomorphic to the closure of the graph of 
the map $\C^2 \dashrightarrow \P^1$ given by $(x,y)\mapsto (x^2:y)$,
over each point in $I(\Phi_n) - I(\Phi_2)$.

We find, 

\begin{theorem} \label{Gamma}
For each $n\geq 2$, there is a regular homeomorphism 
$h_n: B_n\to \Gamma_n$ 
such that $\pi_n \circ h_n = p_n$.  
\end{theorem}

\noindent
In particular, the homeomorphism $h_n$ restricts to the identity on 
the dense open subset $M_2$.
Note that $B_n$ has a regular double point over each point of 
$I(\Phi_n) - I(\Phi_2)$, and therefore, $\Gamma_n$ is singular
for all $n\geq 3$.   
It would be interesting to know if $\Gamma_n$ and
$B_n$ are in fact isomorphic.  
This would follow, for example, if $\Gamma_n$ were known 
to be normal.

Recall that $\hat{M}_2$ is defined to be the inverse limit
of the system $\pi_n: \Gamma_n\to \Gamma_{n-1}$.  
The following immediate 
corollary to Theorem \ref{Gamma}
implies that  the boundary of $M_2$ in $\hat{M}_2$ 
looks (topologically) like the drawing of Figure 1.

\begin{cor}  \label{Bhat}
The inverse limit space $\hat{M}_2$ is naturally homeomorphic
to the inverse limit $\hat{B}$ of the system $p_n: B_n \to B_{n-1}$.
\end{cor}

We are now able to prove Theorem \ref{finite} which states that
any sequence in $\hat{M}_2\subset \prod_{n=1}^\infty \Mbar_{2^n}$
is determined by finitely many entries.

\bigskip\noindent{\bf Proof of Theorem \ref{finite}.}
Let $x = (x_1, x_2, x_2, \ldots)$ denote a sequence in 
$\hat{M}_2$.  Suppose first that $x_1\not\in I(\Phi_n)$ for 
any $n\geq 2$.  Then every iterate map $\Phi_n$ is regular
at $x_1$ so that it has well-defined iterates $x_n\in \Mbar_{2^n}$
for all $n$.  Consequently the sequence $x$ is determined by the
single entry $x_1$.  

Now suppose that $x_1\in\Mbar_2$ is in $I(\Phi_n)$ for some $n\geq 2$,
and let $N$ be the minimal such $n$.  The claim is that 
$x$ is determined by $(x_1, \ldots, x_N)$.  

Indeed, let $y_N$ be
the point in $B_N$ identified with $(x_1, \ldots, x_N)\in \Gamma_N$ via 
the homeomorphism of Theorem \ref{Gamma}.  
By the definition of $B_n$ for each $n\geq N$, the composition of 
projections $p_{N +1} \circ 
\cdots \circ p_n: B_n \to B_N$ is an isomorphism near 
$y_N$.  Consequently, for every $n\geq N$, there is a 
unique point $y_n\in B_n$ associated to $y_N$ such that 
$p_{n+1}(y_{n+1}) = y_n$.  Using again 
Theorem \ref{Gamma}, we find that $(x_1, \dots, x_n)$ is
determined by $(x_1, \dots, x_N)$ for every $n\geq N$, proving the claim.  
\qed

\bigskip\noindent {\bf Proof of Theorem \ref{Gamma}.}
For each $a\in\Chat$, 
let $\Lambda_a\in\Ratbar_2$ be defined
by (\ref{lambda}).  Recall from Theorem \ref{indet} that the indeterminacy
locus of $\Phi_n$ in $\Mbar_2$ is the finite set 
 $$I(\Phi_n) =  \{[\Lambda_\zeta]: \zeta\not=1 \mbox{ and }
\zeta^q=1 \mbox{ for some } q\leq n \}.$$  
For each $q\leq n$ and primitive $q$-th root of unity $\zeta$, 
There exists a regular 
homeomorphism from $\P^1$ to the 
the fiber of the composition $\pi_2\circ \pi_3\circ \cdots \pi_n:
\Gamma_n\to \Mbar_2$ over 
$[\Lambda_\zeta]$:  in local coordinates, the map is given by 
$$\tau^2 \mapsto ([\Lambda_\zeta], [\Lambda_\zeta^2], \ldots, 
[\Lambda_\zeta^{q-1}], [F_{q,\tau,q}], \ldots, [F_{q,\tau,n}])\in\Gamma_n,$$ 
where $\tau$ is a square
root of $\tau^2\in\C$, and 
 $$\infty\mapsto ([\Lambda_\zeta], [\Lambda_\zeta^2], \ldots, 
[\Lambda_\zeta^{q-1}],[P_{q,q}], \dots, [P_{q,n}])\in\Gamma_n.$$  
That this map is well-defined follows from Proposition \ref{limtau}. 
Injectivity follows from Lemma \ref{nonconstant}, and continuity 
from Lemma \ref{continuous}.  The fiber
is not necessarily isomorphic to $\P^1$, as it might have a singularity
at the image of $\infty$.  

To compare $\Gamma_n$ with $B_n$, consider the composition
of projections  $q_n = p_2\circ p_3 \circ \cdots \circ p_n : B_n\to \Mbar_2$.
Observe that
the exceptional fiber of $q_n:B_n\to\Mbar_2$ over $[\Lambda_{-1}]$ is a $\P^1$, 
corresponding
to the family of lines in $\Mbar_2$ passing through the point $[\Lambda_{-1}]$, 
and the exceptional fibers over all other points in $I(\Phi_n)$ are each a $\P^1$ in 
correspondence with a family of conics in $\Mbar_2$ 
passing through $[\Lambda_\zeta]$ which are tangent to $\del M_2$.  

Let $\zeta$ be a primitive $q$-th root of unity, $q\geq 2$.  Suppose 
that $\Delta:\D\hookrightarrow \Mbar_2$ is a
holomorphic disk in $\Mbar_2$ such that 
$\Delta(0) = [\Lambda_\zeta]$.  From Theorem \ref{tau2},
we see that the limiting value
of the iterates, $\Phi_n(\Delta(t))$ as $t\to 0$, depends precisely on the limiting
value, 
  $$\lim_{t\to 0} \frac{(\alpha(t)^q - 1)^2}{\eps(t)} = \tau^2(\Delta)\in\Chat.$$
It suffices to put the exceptional fiber of $B_n$, as a family of lines or conics
through $[\Lambda_\zeta]$, in 
correspondence with the parameter $\tau^2$.  Indeed, we will consider
limits of the iterate map along a given line or conic passing through 
$[\Lambda_\zeta]$ 
and compute the corresponding value of $\tau^2$.  
We will treat the cases $q=2$ and $q>2$ separately.

\medskip\noindent {\bf Coordinates on $\Mbar_2$.}
Choose coordinates $(x_1:x_2:x_3)$ on $\Mbar_2\iso\P^2$ so 
that when $x_3=1$, we have Milnor's coordinates 
$x_1=\sigma_1$ and $x_2=\sigma_2$ on $M_2\iso\C^2$.
In these coordinates, the boundary of $M_2$ is parameterized
by $[\Lambda_a] = (1:a+ 1/a: 0)$ for $a\in\Chat$.  

\medskip\noindent {\bf Case $q=2$.}
Let $\zeta = -1$, so that $[\Lambda_\zeta] = (1:-2:0) \in \Mbar_2$.  
The exceptional fiber of $B_n$ over $[\Lambda_\zeta]$ is the fiber of the projectivized
tangent bundle of $\P^2$ at $[\Lambda_\zeta]$, and can be identified with 
the family of lines 
  $$L_{(a:b)} = \{ (x_1:x_2:x_3): 2b x_1 + b x_2 - a x_3 = 0\},$$
for $(a:b)\in\P^1$.  Each 
line can be parametrized near $[\Lambda_\zeta]$ by 
  $$t \mapsto (1: -2 + a t:  bt)  =: f_t \in\Mbar_2,$$
for $t\in\D$, so that $f_0 = [\Lambda_\zeta]$.  

First suppose that $b=0$ and set $a=1$.  Then $f_t = (1:-2+t:0)$
is in $\del M_2$ for all $t\in\D$.  By definition, 
this direction of approach corresponds
to the parameter $\tau^2 = \infty$.  

Now assume that $b\not=0$.  
For each $t\not=0$, we can write,
  $$f_t = (\sigma_1(t): \sigma_2(t):1) = (1/bt: (-2+at)/bt: 1),$$
so that the fixed point multipliers, $\alpha=\alpha(a,b,t)$, $\beta=\beta(a,b,t)$, and 
$\gamma=\gamma(a,b,t)$, are the three roots of the equation, 
  $$bt x^3 -  x^2 + (at -2) x + 2bt -1 = 0.$$
By construction, there are exactly
two solutions which approach $\zeta = -1$ as $t\to 0$ and one solution
tending to $\infty$.  Label them so that $\gamma\to\infty$ and 
$\alpha + \beta \to -2$. 
While we can't label $\alpha$ and $\beta$ individually for all $t\in\D$, 
we aim to 
find an expression for $(\alpha^2 - 1)^2$ as a function of $t$, 
by using the relations,  
\begin{itemize}
\item[(i)] $\alpha + \beta + \gamma = \frac{1}{bt}$, 
\item[(ii)] $\alpha\beta + \alpha\gamma + \beta\gamma
= \frac{a}{b} - \frac{2}{bt}$, and 
\item[(iii)] $\alpha\beta\gamma = \frac{1}{bt} - 2$.
\end{itemize}
It follows from the (i) that the meromorphic function $\gamma(t)$ can
be expressed as 
  $$\gamma(t) = \frac{1}{bt} + 2 + c_1 t + O(t^2),$$
for some constant $c_1$, and 
  $$(\alpha + \beta) (t) = -2 - c_1 t + O(t^2).$$
Expressions (ii) and (iii) allow us to compute all of the coefficients 
in the expansion of $\gamma$ in terms of $a$ and $b$, but we 
need only $c_1$.  

Expressing $\gamma^{-1}$ as a power series in 
$t$, we obtain
  $$\gamma^{-1} = bt(1 + 2bt + O(t^2))^{-1} = bt(1 - 2bt + O(t^2)),$$
and therefore, (iii) implies that 
  $$\alpha\beta = -\gamma^{-1}(2 - 1/bt) = 1 - 4bt + O(t^2).$$
On the other hand, it follows from (ii) that 
\begin{eqnarray*}
  \alpha\beta &=& -\gamma(\alpha + \beta) + \frac{a}{b} - \frac{2}{bt} \\
&=& \frac{a}{b} + \frac{c_1}{b} + 4 + O(t),
\end{eqnarray*}
and we can solve for $c_1$ to obtain $c_1 = -3b - a$.  Therefore, we can write,
  $$\alpha + \beta = -2 + (3b + a) t + O(t^2).$$
Let us note that $\alpha$ and $\beta$ are the roots of the equation 
  $$x^2 - (\alpha+\beta) x + \alpha\beta = x^2 - (-2 + (3b+a)t + O(t^2)) x 
+ (1 - 4bt + O(t^2)) = 0.$$
According to the quadratic formula, $\alpha$ has the form 
 $$\alpha = -1 + \frac12 (a+3b) t  + O(t^2) \pm \frac12\sqrt{ 4(b-a)t + O(t^2)},$$
and therefore 
  $$(\alpha^2 - 1)^2 = 4(b-a)t + o(t).$$
Finally, we are able to compute 
  $$\tau^2 = \lim_{t\to 0} \frac{(\alpha(t)^2 - 1)^2}{1 - \alpha(t)\beta(t)}
= \frac{b-a}{b}.$$
Consequently, the parameter $\tau^2\in\Chat$ is in one-to-one 
correspondence with the family of lines $L_{(a:b)}$ for $(a:b)\in\P^1$.
This completes the case of $q=2$.

\medskip\noindent {\bf Case $q>2$.}  Let $\zeta$ be a primitive
$q$-th root of unity.  The exceptional fiber of $B_n$ over $[\Lambda_\zeta]
= (1: \zeta + 1/\zeta: 0)$ can be identified with the family of conics,
  $$C_{(a^2:b^2)} = \{ (x_1:x_2:x_3): a^2 x_1x_3 - 
b^2(x_2 - (\zeta + 1/\zeta)x_1)^2 = 0 \},$$
parameterized in a two-to-one fashion by $(a:b)\in\P^1$, each tangent 
at $[\Lambda_\zeta]$ to the boundary $\del M_2$.    
The curve $C_{(a^2:b^2)}$ can be parameterized near $[\Lambda_\zeta]$ by 
  $$t\mapsto f_t = (1: (\zeta + 1/\zeta) + at: b^2t^2)\in\Mbar_2.$$

First suppose that $b=0$ and set $a=1$.  Then $f_t = 
(1: \zeta + 1/\zeta + t: 0)$ is in $\del M_2$ for all $t\in\D$.  
This direction of approach
corresponds to the parameter $\tau^2 = \infty$.  

Now assume that $b\not=0$, so that for each $t\not=0$, we can write, 
  $$f_t = (\sigma_1(t):\sigma_2(t):1) = (1/b^2t^2: (at + \zeta + 1/\zeta)/b^2t^2: 1).$$
The fixed point multipliers, $\alpha = \alpha(a,b,t)$, 
$\beta=\beta(a,b,t)$, and 
$\gamma=\gamma(a,b,t)$, are the three roots of the equation,
 $$b^2t^2 x^3 - x^2 + (at + \zeta + 1/\zeta) x  + 2b^2t^2 - 1 = 0.$$
To compute the value of $\tau^2$, we use the same analysis as in 
the case of $q=2$.  The difference is that all three of the multipliers
can be labelled as meromorphic functions of $t$ such that 
$\alpha \to \zeta$, $\beta \to 1/\zeta$, and $\gamma \to \infty$ 
as $t\to 0$.  Computing the first few terms in the power series for
$\alpha$ and $\beta$  leads to 
  $$\tau^2 = \lim_{t\to 0} \frac{(\alpha(t)^2 - 1)^2}{1 - \alpha(t)\beta(t)}
= \frac{-q^2a^2\zeta^3}{b^2(\zeta^2-1)(\zeta -1)^2}.$$
Therefore, the parameter $\tau^2\in\Chat$ is in one-to-one correspondence
with the family of conics $C_{(a^2:b^2)}$, and this completes the 
proof of Theorem \ref{Gamma}.  
\qed


\bigskip
\section{The space of barycentered measures}
\label{barycenter}

In this section, we study the space of barycentered probablity 
measures on the Riemann sphere and its quotient by the 
group of rotations $SO(3)$.  

Let $M^1(\Chat)$ denote the space of probability measures 
on the Riemann sphere, with the weak-$*$ topology.  
Identify $\Chat$ with the unit sphere in $\R^3$ by stereographic
projection, and let the unit ball in $\R^3$ be taken as a model
for hyperbolic space $\Hyp^3$.  
Recall that the group of M\"obius transformations $\Aut\Chat\iso
\PSL_2\C$ is also the group of orientation preserving isometries
of $\Hyp^3$.

The Euclidean center of mass of a probability measure
$\mu$ on $S^2$ is given by 
  $$E(\mu) =\int_{S^2} \zeta \,d\mu(\zeta).$$ 
Given $\mu\in M^1(\Chat)$ such that $\mu(\{z\}) < 1/2$
for all $z\in\Chat$, the conformal barycenter $C(\mu)\in\Hyp^3$ 
is uniquely determined by the following two properties
\cite{Douady:Earle}:  
\begin{itemize}
\item[1.]  $C(\mu) = 0$ in $\R^3$ if and only if $E(\mu) = 0$, and
\item[2.]  $C(A_*\mu) = A(C(\mu))$ for all $A\in\Aut\Chat$.
\end{itemize}
The barycenter is a continuous function on the space of probability 
measures such that $\mu(\{z\})< 1/2$ for all $z\in\Chat$, and it is 
undefined if $\mu$ has an atom of mass $\geq 1/2$.  
A measure $\mu$ is said to be {\bf barycentered} if  $C(\mu) = 0$.  

  
Let $BCM\subset M^1(\Chat)$ denote the subspace of barycentered measures.  
It is invariant under the action of 
the compact group of rotations $SO(3)\subset \PSL_2\C$.  Let $\overline{BCM}$
denote the closure of $BCM$ in $M^1(\Chat)$.  We will consider the quotient
topological spaces $BCM/SO(3) \subset \overline{BCM}/SO(3)$.  

\begin{theorem} \label{BCM}
The quotient space  $BCM/SO(3)$ is a locally compact
Hausdorff topological space, with the topology induced by the weak 
topology on $M^1(\Chat)$.  The quotient space $\overline{BCM}/SO(3)$ is
the one-point compactification of $BCM/O(3)$.  
\end{theorem}

This first lemma will imply that the quotient of $BCM$ is Hausdorff.  

\begin{lemma}  \label{Hausdorff}
Suppose $\{\mu_k\}$ and $\{\nu_k\}$ are sequences in $BCM$ 
such that $\nu_k = g_{k*} \mu_k$ for a sequence of automorphisms
$g_k\in SO(3)$.  If $\mu_k \to \mu$ and $\nu_k \to \nu$ weakly, then
$\nu = g_* \mu$ for some $g\in SO(3)$.  
\end{lemma}

\proof
Let $\phi$ be a continuous function on $\Chat$.  Passing to a subsequence if 
necessary, we can assume that $g_k\to g$ in $SO(3)$.   
Then $\phi\circ g_k$ converges uniformly to $\phi\circ g$.  Therefore,
the quantity 
  $$\left|\int (\phi\circ g_k) \mu - \int (\phi\circ g) \mu\right|$$
can be made as small as desired for sufficiently large $k$
uniformly over all probability measures $\mu$.   
This estimate together with 
weak convergence of $\mu_k \to \mu$ shows that  
  $$\left|\int (\phi\circ g_k) \mu_k - \int (\phi\circ g) \mu\right|$$
can be made arbitrarily small as $k\to\infty$.  Since this holds
for every $\phi$, we conclude that $g_{k*}\mu_k \to g_*\mu$
weakly.  On the other hand, $\nu_k = g_{k*}\mu_k \to \nu$
weakly, so $\nu = g_*\mu$.  
\qed

\bigskip
With the identification of $\Chat$ and the unit sphere in $\R^3$, 
note that the antipode of a point $a\in\Chat$ is $-1/\bar{a}$.  
The following lemma shows what happens to an unbounded sequence
in $BCM$.  

\begin{lemma}  \label{1/2}
Let $\mu_k$ be a sequence of barycentered measures such that 
$\mu_k \to \mu$ weakly.  Then either $\mu$ is barycentered, or
  $$\mu = \frac12 \delta_a + \frac12 \delta_{-1/\bar{a}}.$$
\end{lemma}

\proof
Suppose first that $\mu(\{z\}) < 1/2$ for all $z\in\Chat$.  Then 
the barycenter of $\mu$ is well-defined.  It follows from
the continuity of the barycenter that $\mu$ is barycentered.  
Now suppose there is a point $a\in\Chat$ such that 
$\mu(\{a\}) \geq 1/2$.  Then by weak convergence, there
exists a sequence $\eps_k \to 0$ and $r_k \to 0$ such that 
  $$\mu_k(B(a,r_k)) \geq \frac{1}{2} - \eps_k.$$

Without loss of generality, we may assume that $a = (1,0,0)$ in $\R^3$.
Then $-1/\bar{a} = (-1,0,0)$.  
It suffices to show that for each fixed $r>0$, there is a sequence 
$\delta_k \to 0$ so that 
  $$\mu_k (B(-1/\bar{a},r)) \geq 1/2 - \delta_k.$$
Then the limiting measure $\mu$ must satisfy 
  $$\mu(\{a\}) = \mu(B(-1/\bar{a},r)) = 1/2,$$
for all $r>0$.  Letting $r\to 0$, we see that $\mu(\{-1/\bar{a}\}) = 1/2$.  

Suppose, upon passing to a subsequence if necessary, that there is 
a $\delta>0$ such that 
  $$\mu_k(B(\hat{a},r)) \leq 1/2 - \delta$$
for all $k$.  Let $\zeta_x$ denote the $x$-coordinate of a vector 
$\zeta\in S^2$.  
Since the measures $\mu_k$ are barycentered, the Euclidean
center of mass of $\mu_k$ is at the origin, and therefore, 
\begin{eqnarray*}
 0 = \int_{\Chat} \zeta_x \mu_k(\zeta) 
&=& \int_{B(a,r_k)} \zeta_x \mu_k(\zeta) 
+ \int_{\Chat - B(a,r_k)\cup B(\hat{a},r)} \zeta_x \mu_k(\zeta)
+ \int_{B(\hat{a},r)} \zeta_x \mu_k(\zeta)  \\
&\geq&  (1/2 - \eps_k)(\cos \pi r_k) + (\delta + \eps_k)(-\cos \pi r)
- (1/2 - \delta)  \\
&=&   \delta(1 - \cos \pi r) + 1/2 (\cos \pi r_k - 1) 
- \eps_k(\cos \pi r + \cos \pi r_k).
\end{eqnarray*}
For sufficiently large $k$, the final line is positive, which
is a contradiction.  
\qed

\bigskip\noindent
{\bf Proof of Theorem \ref{BCM}}.  That the quotient 
$BCM/SO(3)$ is Hausdorff follows from Lemma \ref{Hausdorff}. 
Local compactness of the metrizable space $BCM$ and the 
fact that the $SO(3)$-orbits are compact implies that 
$BCM/SO(3)$ is also locally compact.  
Suppose $\mu_k$ is an unbounded sequence in $BCM$ such 
that $\mu_k \to \nu$ weakly in $M^1(\Chat)$.  By Lemma \ref{1/2}, 
$\nu = \frac12 \delta_a + \frac12 \delta_{-1/\bar{a}}$ for an antipodal pair
$(a,-1/\bar{a})$.  Under the action of $SO(3)$ on $\M^1$, all such 
$\nu$ are equivalent.
\qed

\bigskip\noindent{\bf The point at infinity.}
We will refer to the point at infinity of $\overline{BCM}/SO(3)$ simply
by $\infty$.  One way to detect if a sequence of probability measures
is converging to $\infty$ in $\overline{BCM}/SO(3)$ is to find 
a sequence of ``separating annuli", annuli of growing modulus such that
half of the measure lies on side of the annulus and half on the other.  
By classical arguments, at least one complementary component 
of the annulus must shrink to a point.  

\begin{lemma}  \label{infinity}
Suppose $\{\mu_k\}$ is a sequence of barycentered probability 
measures on $\Chat$ such that $\mu_k \to \nu$ weakly, 
and $A_k$ a sequence of round annuli such that 
\begin{itemize}
\item[(1)] $\mod A_k \to \infty$ as $k\to\infty$, and 
\item[(2)] there is a sequence $\eps_k \to 0$ such that 
$\mu_k(D_k) \geq 1/2 - \eps_k$ for each of the 
complementary disks $D_k$ of $A_k$.
\end{itemize}
Then $\nu = \infty$ in $\overline{BCM}/SO(3)$.
\end{lemma}

\proof
By condition (1), some subsequence of
the closed disks $D_k$ is converging to a point in the Hausdorff
topology, say $D_k \to \{a\}$.  Then for any $r>0$, we have 
  $$\nu(B(a,r)) \geq \lim_{k\to\infty} \mu_k(D_k) \geq 1/2.$$
This holds for all $r>0$, and therefore $\nu(\{a\}) \geq 1/2$,
so by Lemma \ref{1/2}, $\mu_k\to \infty$ in $\overline{BCM}/SO(3)$.
\qed

\bigskip

\section{The homeomorphism $\hat{M}_2\to X_2$}
\label{maintheorem}

In this section we give the proof of Theorem \ref{homeo} which 
states that the compactification by barycentered measures and
the inverse limit space which resolves the iterate maps are the same
in degree 2.  
We begin by studying the boundary behavior of the continuous map 
 $$M: M_2 \to BCM/SO(3),$$ 
which sends a conjugacy class $[f]$ to the maximal measure 
$\mu_f$ of a barycentered representative.  

Recall that the boundary of $M_2$ in $\Mbar_2\iso\P^2$ 
is parametrized
by the family of conjugacy classes $[\Lambda_a] = [\Lambda_{1/a}]$
for $a\in\Chat$, where $\Lambda_a\in\Ratbar_2$ is defined
in (\ref{lambda}).  

\medskip
\begin{prop} \label{cont1}
Suppose $[f_k]$ is a sequence in $M_2$ such that 
$[f_k] \to [\Lambda_a]$ in $\Mbar_2$ as $k\to\infty$, 
where $a$ is not a root of unity.  Then $M([f_k]) \to \infty$
in $\overline{BCM}/SO(3)$ as $k\to\infty$.  
\end{prop}

\proof
First suppose that $a$ is neither 0 nor $\infty$.  
There is a sequence
of representatives $f_k\in\Rat_2$ such that $f_k\to \Lambda_a$ in 
$\Ratbar_2$, so that $\mu_{f_k} \to \mu_{\Lambda_a}$ weakly by 
Theorem \ref{limitmu}.  Let $\mu_a = \mu_{\Lambda_a}$.  
By definition,  
  $$\mu_a = \frac12 \sum_{n=0}^\infty \frac{1}{2^n} \delta_{1/a^n}.$$
Since $a$ is not a root of unity, we have $\mu_a(\{z\})\leq
1/2$ for all $z\in\Chat$, and there exists a point $p\in\Chat$ such 
that $\mu_a(\{p\}) = 1/2$. 

Fix $\eps>0$ and choose $r = r(\eps) > 0$ so that 
  $$\mu_a(\Chat - \overline{B(p,r)}) \geq 1/2 - 2\eps.$$
By weak convergence of the measures $\mu_k \to \mu_a$,
there exists an integer $N(\eps)$ such that 
  $$\mu_k(\Chat - \overline{B(p,r)}) \geq 1/2 - \eps,$$
and
  $$\mu_k(B(p,r^2)) \geq 1/2 - \eps$$
for all $k\geq N(\eps)$.  
We can assume that $r(\eps) \to 0$ and $N(\eps) \to \infty$
as $\eps \to 0$.  

Rephrasing, given $k$, we can let $\eps_k$ be
the smallest $\eps$ such that $k\geq N(\eps)$, and set
$r_k = r(\eps_k)$.  Then as $k\to\infty$, we have 
$r_k\to 0$.  Consequently, the 
annulus $B(p,r_k) - \overline{B(p,r_k^2)}$ has
$\mu_k$-measure $< 2 \eps_k$ and modulus $\to \infty$
as $k\to \infty$.

Now suppose that $g_k\in \Aut\Chat$ is chosen so that $g_{k*}\mu_k$ is
barycentered.  Let $A_k = g_k(B(p,r_k) - \overline{B(p,r_k^2)})$,
so that $\mod A_k \to \infty$.  
If $\nu$ is any subsequential limit of the measures  $g_{k*}\mu_k$,
then Lemma \ref{infinity} implies that $\nu = \infty$ in 
$\overline{BCM}/SO(3)$.  

Finally, suppose that $a=0$ or $a=\infty$.  
Then the probability measure associated to $\Lambda_a\in\Ratbar_2$
is
 $$\mu_{\Lambda_a} = \frac12 \delta_1 + \frac12 \delta_{-1},$$
and $\mu_k \to \mu_{\Lambda_a}$ weakly.  
Therefore, for any $r>0$, any annulus of the form $\Chat - (B(1,r)\cup
B(-1,r))$ will have $\mu_k$-measure tending to 0 as $k\to \infty$.  
By Lemma \ref{infinity}, the barycentered representatives of $[f_k]$ must
tend to infinity in $\overline{BCM}/SO(3)$.
\qed

\bigskip
Let $\zeta$ be a primitive $q$-th root of unity for some $q\geq 2$.
Recall that by Theorem \ref{tau2}, the limiting value of 
the iterate map $\Phi_n$ near $[\Lambda_\zeta]\in\del M_2$ depends on the direction 
of approach, and the $n$-th iterate can be computed in terms 
of a parameter $\tau^2\in\Chat$.  
The limiting barycentered measure depends 
on the limit of the iterates.  The key observation is that the measures 
associated to the 
stable limits of the $q$-th iterates (in $\Ratbar_{2^q}$) 
have well-defined barycenters.  Recall the definitions of $F_{q,\tau,q}$
and $P_{q,q}\in\Ratbar_{2^q}$ 
from Propositions \ref{Fiterates} and \ref{paraboliclimit}.  

\medskip
\begin{prop} \label{cont2}
Fix $q\geq 2$ and $\zeta$ a primitive $q$-th root of unity.  
Let $[f_k]$ be a sequence in $M_2$ converging to $[\Lambda_\zeta]$
in $\Mbar_2$ as $k\to\infty$ such that the $q$-th iterates $[f_k^q]$ converge 
in $\Mbar_{2^q}$ to either (i) $[P_{q,q}]$ or 
(ii) $[F_{q,\tau,q}]$ for some $\tau\in\C$.  
Then in case (i), $M([f_k]) \to \infty$ in $\overline{BCM}/SO(3)$ as $k\to \infty$.
In case (ii), 
$\lim_{k\to \infty} M([f_k])$ is equivalent (in $\PSL_2\C$) to the measure $\mu_{F_{q,\tau,q}}$.
\end{prop}

\proof
In case (i), there exist representatives $f_k\in\Rat_2$ such that the
iterates $f_k^q$ converge in $\Ratbar_{2^q}$ to $P_{q,q}$ as $k\to\infty$
Since $P_{q,q} \not\in I(2^q)$, Theorem 
\ref{limitmu} implies that the measures $\mu_{f_k}$ converge weakly to 
  $$\mu_{P_{q,q}} = \frac{1}{2^q} \sum_{j=0}^\infty \frac{2^{q-1}}{2^{qj}}
\delta_{-j}  +  \frac{1}{2^q} \sum \frac{2^{q-1}-1}{2^{qj}} \delta_\infty 
= \frac12 \sum_{j=0}^\infty \frac{1}{2^{qj}} \delta_{-j} 
+ \frac{2^{q-1}-1}{2^q-1} \delta_\infty$$
as $k\to\infty$.  Notice that $\mu_{P_{q,q}} (\{0\}) = 1/2$. 

Therefore, there exists a family of annuli  $A_k$ separating
$z=0$ from the other points in the support of $\mu_{P_{q,q}}$,  
such that $\mod A_k \to \infty$ as $k\to\infty$ and  the $\mu_{f_k}$-measure
of each complementary component of $A_k$ tends to $1/2$ as
$k\to\infty$.  Choose a sequence $g_k\in\Aut\Chat$ so that each 
$g_{k*} \mu_{f_k}$ is barycentered.  Then $g_{k*} \mu_{f_k}$ and
the annuli $g_k(A_k)$ satisfy the hypotheses of Lemma \ref{infinity}
which shows that $g_{k*}\mu_{f_k} \to \infty$ in $\overline{BCM}/SO(3)$.
 
In case (ii), there exist representatives $f_k\in\Rat_2$ such that
the iterates $f_k^q$ converge in $\Ratbar_{2^q}$ to $F_{q,\tau,q}$
as $k\to\infty$.  Since $F_{q,\tau,q}$ is not in $I(2^q)$, Theorem
\ref{limitmu} implies that the measures $\mu_{f_k}$ converge weakly 
to $\mu_{F_{q,\tau,q}}$ as $k\to\infty$.  Note that $\mu_{F_{q,\tau,q}} (\{z\}) 
< 1/2$ for all $z\in \Chat$ by Lemma \ref{Fmeasures}.  
Therefore, $\mu_{F_{q,\tau,q}}$ 
has a well-defined barycenter.  By continuity of the barycenter, 
we can choose a sequence $g_k\in\Aut(\Chat)$ for 
such that $g_{k*} \mu_{f_k}$ is barycentered for all $k$ and 
$g_k\to g \in\Aut(\Chat)$ with $g_* \mu_{F_{q,\tau,q}}$ barycentered.
Therefore, $M([f_k])$ converges 
in $\overline{BCM}/SO(3)$ to a measure which is equivalent 
(in $\Aut(\Chat)$) to 
$\mu_{F_{q,\tau,q}}$ as $k\to\infty$.  
\qed

\bigskip\noindent
{\bf Proof of Theorem \ref{homeo}}.  We aim to show that 
the embedding of $M_2$ into $M_2\times BCM/SO(3)$ via
the graph of $M([f]) = \mu_f$ (for a barycentered representative)
extends to a homeomorphism  
 $$h: \hat{M}_2 \to X_2 \subset \Mbar_2\times \overline{BCM}/SO(3).$$
Since $\hat{M}_2$ is compact and $X_2$ is Hausdorff, 
it suffices to show that $h$ is 
continuous and bijective.  Furthermore, since its image contains a dense
open subset, namely $M_2$ itself, it suffices to show only continuity 
and injectivity.  

Let $x = (x_1, x_2, \dots)$ be a boundary point of $M_2$ 
in $\hat{M}_2\subset\prod_{n=1}^\infty \Mbar_{2^n}$.  Suppose 
first that $x_1 = [\Lambda_a]$ where $a$ is not a primitive $q$-th
root of unity for any $q\geq 2$.  Then by Theorem \ref{indet}, 
$x_n = [\Lambda_\zeta^n]$ for every $n\geq 2$.  That is to say,
there is a unique point in $\hat{M}_2$ which projects to $x_1$
in $\Mbar_2$.  By Proposition \ref{cont1}, $h$ extends continuously
to $x$ with $h(x) = ([\Lambda_a], \infty) \in 
\Mbar_2\times\overline{BCM}/SO(3)$, and we see that $h(x)\in X_2$ has
$x$ as its unique preimage.  

Now suppose that $x_1\in I(\Phi_n)$ for some $n\geq 2$.  Then 
by Theorem \ref{indet}, 
there is a $q\leq n$ and a primitive $q$-th root of unity $\zeta$ such 
that $x_1 = [\Lambda_\zeta]$.  From Proposition \ref{limtau}, 
the $q$-th entry $x_q$ in the sequence $x$ 
must be either (i) $[P_{q,q}]$, or (ii) of the form $[F_{q,\tau,q}]$ 
for some $\tau\in\C$, and all further entries $x_n$ for $n> q$ are 
determined by $x_q$.  That is to say, there is a unique point in $\hat{M}_2$
which projects to $(x_1,\ldots, x_q)$ under the projection to the first
$q$ factors.  
Let $[f_k]$ be any sequence in $M_2$ such that 
$[f_k] \to x$ in $\hat{M}_2$ as $k\to \infty$.  By Proposition \ref{cont2}, 
$h$ extends continuously to $x$ with $h(x) = ([\Lambda_\zeta], \infty)$
in case (i) and $h(x) = ([\Lambda_\zeta], g_*\mu_{F_{q,\tau,q}})$ for 
some $g\in\Aut\Chat$ in case (ii).  Furthermore, we see that $h(x)\in X_2$
has $x$ as its unique preimage by Lemma \ref{nonconstant}.  
\qed

\bigskip
\section{Higher degrees} 
\label{higher}

In this section, we show that the inverse limit space 
$\hat{M}_d$ and the compactification by barycentered measures
$X_d$ are not  homeomorphic for $d\geq 5$.  The
examples come from \cite[\S5]{D:measures}.  We also 
prove the following theorem.

\begin{theorem}  \label{discontinuity}
The iterate map $\Phi_n : M_d \to M_{d^n}$ does 
not extend continuously to $\Mbar_d$ for any $d\geq 2$
and $n\geq 2$.  
\end{theorem}

The examples used to show discontinuity 
at the GIT boundary of $M_d$ are also from \cite[\S5]{D:measures}.
In fact, they are a generalization of Epstein's examples
for the proof of Proposition \ref{Prop2} 
(\cite[Prop 2]{Epstein:bounded}).  The idea is to find
unbounded families in $\Rat_d$ such that the  
critical points are at $2d-2$ given points for every map in the 
family, and such that these critical points are distinct
from the holes which develop in the limit.

\bigskip\noindent {\bf The second iterate.}
We give first a complete proof of Theorem \ref{discontinuity}
for the case of the second iterate, $n=2$.  

Fix $d\geq 2$.  
Let $P(z,w)$ be a homogeneous polynomial of degree
$d-1$ with distinct roots in $\P^1$ which is monic as 
a polynomial in $z$ and such that 
$P(0,1)\not=0$ and $P(1,0)\not=0$.  Let 
   $$g=(wP(z,w): 0) \in\Ratbar_d.$$  
Then $g$ has
a hole at $\infty$ of depth 1 and holes of depth 1 at 
each of the roots of $P$.  The lower degree map $\phi_g$
is the constant $\infty$ map, so we see that $g\in I(d)$.  
The point $g$ is stable for all $d\geq 4$,
semistable for $d=3$, and unstable for $d=2$.  

Consider the family of rational
maps given by 
  $$g_{a,t} = (at z^d + wP(z,w): tz^d),$$
for $a\in\C$ and $t\in\D^*$.  This family converges 
to $g$  in $\Ratbar_d$ as $t\to 0$ for every $a\in\C$.  
In \cite[\S5]{D:measures}, it
was computed that the second iterates of 
the family $g_{a,t}$ converge as $t\to 0$ to 
  $$f_a = (w^{d-1} P^{d-1} (awP + z^d) : w^dP^d)\in\Ratbar_{d^2}.$$
In the notation $f_a = H_a\phi_a$, 
we have $H_a = w^{d-1}P^{d-1}$ for all $a\in\C$ and
  $$\phi_a = (awP(z,w) + z^d: wP(z,w))\in\Rat_d.$$
The point $f_a$ has holes at $\infty$ and the roots 
of $P$, each of depth $d-1$.  Note that $f_a$ is {\em not}
in $I(d^2)$.

\medskip
\begin{lemma} \label{gfamily}
For each $d\geq 2$, the conjugacy classes $[g_{a,t}]$ 
converge in $\Mbar_d$ as $t\to 0$ to a boundary point
independent of $a\in\C$.  
\end{lemma}

\proof
For each $d>2$, the point $g\in\Ratbar_d$ is stable
or semistable and therefore determines a unique point
$[g] \in\Mbar_d$.  Convergence of $g_{a,t}$ to $g$ in
$\Ratbar_d$ as $t\to 0$ implies that $[g_{a,t}] \to [g]$ 
in $\Mbar_d$ as $t\to 0$ for every $a\in\C$.  

For $d=2$, however, the point $g$ is unstable, so
it does not represent a point in $\Mbar_2$.  
Conjugating the family $g_{a,t}$ by 
$A_t(z) = t^{1/2} z$, we obtain new representatives, 
  $$A_tg_{a,t}A_t^{-1} = (at^{1/2} z^2 + wP(z, t^{1/2}w): z^2),$$
which converge as $t\to 0$ to 
  $$h = (zw: z^2).$$
This point $h$ is stable, since it agrees with the lower 
degree map $\phi_h(z) = 1/z$ away from a 
hole of depth 1 at $z=0$ which is not fixed by $\phi_h$.
Therefore, $[g_{a,t}] \to [h]$ in $\Mbar_2$ 
as $t\to 0$ for every $a\in\C$.  
\qed

\medskip
\begin{lemma} \label{f_a}
For each $d\geq 2$, 
\begin{itemize}
\item[(i)] $f_a\in\Ratbar_{d^2}$ is stable for all $a\in\C$, and
\item[(ii)] the map $\C \to \Mbar_{d^2}$ given by $a\mapsto [f_a]$ is 
            non-constant.  
\end{itemize}
\end{lemma}

\proof
From the definition of $f_a$ we see that each hole has depth 
$d-1$ and $d-1 = (d^2-1)/(d+1) < (d^2-1)/2$ for all $d\geq 2$.
This proves (i).  
Consequently, $[f_a] = [f_b]$ in $\Mbar_{d^2}$ if and only if 
$f_a$ and $f_b$ are conjugate by an element of $\PSL_2\C$.

Now suppose that $d>2$.  Any conjugacy between $f_a$ and $f_b$ 
for $a\not=b$ must preserve the holes at $\infty$ and
the roots of $P$ and conjugate $\phi_a$ to $\phi_b$.  
There are at least three
holes since the degree of $P$ is $d-1$ with distinct
roots.   
On the other hand, the finite fixed points of
$\phi_a$ are the $d-1$ solutions to $z^d - zP(z,1) + aP(z,1) =0$,
and so the set of fixed points varies with $a\in\C$.
The cross-ratio of three of the holes with a moving fixed
point must then vary with $a\in\C$, and this proves that not
all $f_a$ are conjugate.  

For $d=2$, note that a conjugacy between $f_a$ and $f_b$ 
must preserve the holes and also send the critical points of 
$\phi_a$ to the critical points of $\phi_b$.  It can be computed 
directly that if $P(z,1) = z-\alpha$, 
the critical points of $\phi_a$ are at $z=0$ and $z=2\alpha$, 
independent of $a\in\C$.  Together with 
the root $\alpha$ of $P$ and the point at $\infty$, there
are four marked points which must be permuted by any conjugacy.  
Moreover, the finite fixed point of 
$\phi_a$ is at $z = a\alpha/(a+\alpha)$, so the cross ratio
of $0$, $\alpha$, $\infty$, and the finite fixed point depends
on $a\in\C$.  We conclude that not all $f_a$ are conjugate, and 
the lemma is proved.  
\qed

\medskip
\begin{cor}
The second iterate map $\Phi_2:M_d \to M_{d^2}$ does
not extend continuously to $\Mbar_d$.  
\end{cor}

\proof
This is immediate from Lemmas \ref{gfamily} and \ref{f_a}.
\qed

\medskip

\bigskip\noindent {\bf Higher iterates of $g_{a,t}$}.  
We are now ready to complete the proof of Theorem
\ref{discontinuity}.  

\begin{lemma} \label{stability}
For each $d\geq 2$, $n\geq 2$, and $a\in\C$, the limit of the 
iterates $(g_{a,t})^n\in\Ratbar_{d^n}$ as $t\to 0$ is stable.  
\end{lemma}

\proof   Let $n$ be an even integer. 
We have seen that the second iterates of $g_{a,t}$ converge
to $f_a\in\Ratbar_{d^2}$ as $t\to 0$.  The $n$-th iterates of 
the family $g_{a,t}$ will converge to $(f_a)^{n/2}$
as $t\to 0$ by the continuity of the $n/2$-th iterate map
at $f_a \not\in I(d^2)$.
In \cite[\S5]{D:measures}, it was computed that 
$\mu_{f_a}(\{\infty\}) = 1/(d+1)$ (or it follows from Lemma \ref{Lemma5}).  
The roots of $P$ 
are simple and are each mapped to $\infty$ by $\phi_a$ 
with multiplicity 1, so Lemma \ref{Lemma5} implies that 
$\mu_{f_a}(\{\alpha\}) = 1/(d+1)$ for each root $\alpha$
of $P$.    Since $\mu_{f_a}$ is
a probability measure, any other point of $\P^1$ 
must have mass $\leq 1/(d+1)$.  By Propositions \ref{evendeg}
and \ref{odddeg},
we see that every iterate of $f_a$
must be GIT stable.  Therefore, all even iterates 
of the family $g_{a,t}$ have a stable limit as $t\to 0$.

Now let $n\geq 3$ be an odd integer. 
Lemma \ref{composition} implies that the 
composition map $\scriptC_{d,d^{n-1}}$ is continuous
at the pair $(g, (f_a)^{(n-1)/2})$.  Consequently, 
the $n$-th iterates of $g_{a,t}$ converge to the point 
$g\circ (f_a)^{(n-1)/2}$ as $t\to 0$.  The iterate formula
of \cite[Lemma 7]{D:measures} (see also \S\ref{ratd}) 
shows that
 $$(f_a)^{(n-1)/2} = \prod_{k=0}^{\frac{n-1}{2}-1} 
\left( (\phi_a^k)^*H_a \right)^{d^{- k - 1 + (n-1)/2}} 
(\phi_a)^{(n-1)/2}.$$
It will be useful to write $(\phi_a)^{(n-1)/2}$ in terms of
its coordinate functions $(\phi^{(n-1)/2}_{az}:\phi^{(n-1)/2}_{aw})$
so that we can compute the composition $g\circ (f_a)^{(n-1)/2}$ .  
Indeed, substituting the coordinate functions for this 
iterate of $f_a$ into the formula for $g$, we obtain 
 $$g\circ (f_a)^{(n-1)/2} = \left( \prod_{k=0}^{\frac{n-1}{2}-1}
\left( (\phi_a^k)^*H_a \right)^{d^{-k + (n-1)/2}} 
\phi^{(n-1)/2}_{aw}
P(\phi^{(n-1)/2}_{az},\phi^{(n-1)/2}_{aw}) : 0 \right),$$
where we have factored out all appearances of $H_a$.  
Notice, in particular, that the expression involving $H_a$ appears
as the $d$-th power of the same expression in $(f_a)^{(n-1)/2}$.
The estimates on $\mu_{f_a}$ given in the previous
paragraph together with 
Lemma \ref{Cor8} imply that the depth of any point for
$(f_a)^k$ is no greater than $d^{2k}/(d+1)$. 
Therefore, the depths of the holes of the composition
$g\circ (f_a)^{(n-1)/2}$
at $\infty$ or at the roots of $P$ will not exceed $d(d^{n-1}/(d+1)) +1$,
where the added 1 comes from each of these holes being 
a simple zero of $\phi^{(n-1)/2}_{aw}$.  Notice that this bound is 
$< d^n/2$ for all even $d$ and $< (d^n-1)/2$ for all odd $d$.
Therefore the holes at $\infty$ and the roots of $P$ do not 
violate the stability criteria.  

For any other point in $\P^1$, we 
know from the above that its depth as a hole of $(f_a)^{(n-1)/2}$
is no more than $d^{n-1}/(d+1)$; its depth as a hole
of this composition cannot  then exceed $d^n/(d+1) + d^{(n-1)/2}$, 
where the second term is the degree of $(\phi_a)^{(n-1)/2}$.
This upper bound on the depth is less than $d^n/2$ except 
when $d=2$ and $n=3$.  In this special case, 
it is easy to check that the point $g\circ f_a \in\Ratbar_8$ 
has holes of depth 3 at $\infty$ and the root of $P$, and a hole
of depth at most 2 at any other point.  
It follows that the composition $g\circ (f_a)^{(n-1)/2}$ is always
stable.  
\qed 

\medskip

\bigskip\noindent
{\bf Proof of Theorem \ref{discontinuity}}.
Let $f_{a,n}\in\Ratbar_{d^n}$ denote 
the limit of the iterates $(g_{a,t})^n$ as $t\to 0$.
By Lemma \ref{stability}, $f_{a,n}$ is stable for all 
$n\geq 2$ and all $a\in\C$, and therefore it determines 
a unique point $[f_{a,n}]$ in $\Mbar_{d^n}$.  
Furthermore, stability
implies that $[f_{a,n}] = [f_{b,n}]$ if and only if 
$f_{a,n}$ and $f_{b,n}$ are conjugate.

From Lemma \ref{gfamily}, the family $[g_{a,t}]$ 
converges in $\Mbar_d$ as $t\to 0$ to a point independent 
of $a\in\C$, for every $d\geq 2$.  Lemma \ref{f_a} implies
that $a\mapsto [f_{a,2}]$ is non-constant.  To 
conclude the proof, we need to show that $a\mapsto [f_{a,n}]
\in\Mbar_{d^n}$ is non-constant for all $n\geq 2$.  

Suppose first that $n$ is even.  Any conjugacy between 
$f_{a,n}$ and $f_{b,n}$ for $a\not=b$ must preserve 
the holes and conjugate $\phi_a^{n/2}$ to $\phi_b^{n/2}$.
In particular, it must preserve the critical points of $\phi_a$,
located at the 
$2d-2$ solutions to $dP(z,1)z^{d-1} - P'(z,1)z^d =0$, 
independent of $a\in\C$.  Together with the hole at $\infty$,
these give at least three marked points to be permuted 
by a conjugacy.  On the other hand, the finite fixed points
of $\phi_a$ are the $d-1$ solutions to $z^d - zP(z,1) + aP(z,1) 
=0$, and the set of these will vary with $a\in\C$.  Consequently,
the cross-ratio of three of the marked points with a moving
fixed point of $\phi_a$ also varies with $a\in\C$, so not 
all $f_{a,n}$ are conjugate.  

Now suppose that $n\geq 3$ is odd.  Any conjugacy between
$f_{a,n}$ and $f_{b,n}$ must preserve holes of the same 
depth.  The formula for $f_{a,n} = g\circ f^{(n-1)/2}_a$,
given in the proof of Lemma \ref{stability}, shows that 
for each $d>2$, there are at least three 
holes at $\infty$ and 
the roots of $P$ which are of the same depth and do not
depend on $a\in\C$.  If $\alpha$ is 
a root of $P$, then the preimages of $\alpha$ by $\phi_a$ 
are also holes of $f_{a,n}$ and do depend on $a\in\C$.  
Therefore, the cross-ratio of $\infty$ with two roots of 
$P$ and a moving preimage of $\alpha$ must vary with $a\in\C$,
and so not all $f_{a,n}$ are conjugate.  For $d=2$, if $P(z,1)
= z-\alpha$, note that the cross-ratio of the holes at $\infty$
and $\alpha$ with the pair of preimages of $\alpha$ by $\phi_a$
is given by 
 $$\chi(a) = \frac{a+\alpha + \sqrt{(a-\alpha)^2 + 4\alpha(a-\alpha)}}
                  {a+\alpha - \sqrt{(a-\alpha)^2 + 4\alpha(a-\alpha)}},$$
which depends on $a\in\C$.  Therefore, not all $f_{a,n}$ are conjugate.
\qed

\bigskip\noindent{\bf The spaces $\hat{M}_d$ and $X_d$.}
We conclude by demonstrating that there cannot exist a continuous
map $X_d \to \hat{M}_d$ which restricts to the identity on $M_d$,
for every $d\geq 5$.  It is likely that there exists a continuous map 
in the opposite direction.  

The following examples are from \cite[\S5 Example 2]{D:measures}.
Fix $d\geq 5$.  Let $P = P(z,w)$ be a homogeneous polynomial of degree 
$d-2$ with distinct roots such that $P(0,1)\not= 0$, $P(1,0)\not=0$, and 
$P$ is monic as a polynomial in $z$.  Let $g = (w^2P(z,w):0)$.  Then $g\in I(d)$
is stable for all $d\geq 6$ and semistable for $d=5$ since the depth 
at $\infty$ is $2 < d/2$.  Therefore $g$ defines a unique point $[g]$ in 
$\Mbar_d$.  

For each $a\in\C$ and $t\in[0,1]$, 
consider the family
  $$h_{a,t} = (atz^d + w^2P(z,w): tz^d) \in \Rat_d.$$
Computing second iterates and taking a limit as $t\to 0$, we obtain,
  $$(h_{a,t})^2 \to h_a := (a w^{2d} P(z,w)^d: w^{2d} P(z,w)^d)\in \Ratbar_{d^2}.$$
Note that $h_a$ is stable for all $a\in\C$ and all $d\geq 5$ since the depth 
at $\infty$ is $2d < d^2/2$.  Therefore each $h_a$ determines a point 
$[h_a]\in\Mbar_{d^2}$ and $[h_a]=[h_b]$ if and only if they lie in the same 
$\PSL_2\C$-orbit.  
Write $h_a = H_a\phi_a$.  Since $P$ has at least 3 distinct roots and 
the constant $\phi_a \equiv a$ depends on $a$, we see that only finitely 
many of the conjugacy classes $[h_a]$ can coincide.  

For each $a\in\C$ such that $P(a,1) \not=0$, the point 
$h_a$ is not in $I(d^2)$.  Therefore, the measures $\mu_{h_{a,t}}$
converge weakly as $t\to 0$ to 
  $$\mu_{h_a} = \frac{2}{d} \delta_\infty + \frac{1}{d} \sum_{P(z,1)=0} \delta_z,$$
a measure which is independent of $a$.  The measure $\mu_{h_a}$ has no
atoms of mass $\geq 1/2$, and so it has a well-defined barycenter.  Let
$\mu = g_*\mu_{h_a}$ be a barycentered measure for some $g\in\PSL_2\C$.
Consequently, for every $a\in\C$, the family $[h_{a,t}]$ converges 
in $X_d$ as $t\to 0$ to the pair $(g, \mu) \in \Mbar_d \times \overline{BCM}/SO(3)$.
On the other hand, the limits of $[h_{a,t}]$ as $t\to 0$ in $\hat{M}_d$ are 
distinct since the second iterates have distinct limits in $\Mbar_{d^2}$.  

  
\bigskip
\section{Final remarks}
\label{final}

The choices made in the definitions of $\hat{M}_d$ and $X_d$
reflect the following desirable properties in degree 2.  A compactification
$X$ of the moduli space of quadratic rational maps should satisfy:
\begin{itemize}
\item  iteration is well-defined on $X$, 
\item  there exists a projection from $X$ to Milnor's $\Mbar_2\iso \P^2$,
\item  there exists a projection from $X$ to the space of barycentered 
measures $\overline{BCM}/SO(3)$.
\end{itemize}
In degree 2, the space $\hat{M}_2 = X_2$ can be described as 
the ``minimal'' compactification
satisfying these  properties.  The second condition ensures that there
is a projection of the boundary of $M_2$ to the moduli space $\Mbar_1$. 
The third condition about the barycentered
measures is chosen to reflect the geometry of the rational maps.  Each 
rational map $f$ (together with its measure of maximal entropy) determines 
a convex surface (up to scale) in $\R^3$, with curvature 
equal to the distribution $4\pi\mu_f$, and invariant under conjugation by 
a M\"obius transformation.  
See, for example, \cite[\S6]{D:measures}.  Given an unbounded family in 
$M_2$ and the associated family of surfaces in $\R^3$, the choice of 
barycentered representatives corresponds to fixing the diameter of the surfaces.  
In this language, it would follow from the propositions of \S\ref{maintheorem}
that the limiting geometry of the boundary points are either compact
convex polyhedra
(with countably many vertices) or degenerate needles (where the curvature
is concentrated at two points).  

In higher degrees, it is not immediately obvious how to formulate the ``desirable" properties of a compactification, as very little is known about parametrizations
of the space $M_d$.  

Finally, much can be said about the dynamical properties of unbounded 
families in $M_2$, particularly when restricted to a given hyperbolic component. 
It would be interesting to understand better the explicit examples given in 
Section \ref{degree2proofs} which appeared first in \cite{Epstein:bounded}.  

\medskip\noindent{\bf Example.}  
Suppose $\{f_t: t\in (0,1]\}$ is a family of quadratic rational maps with 
$f_1(z) = z^2-1$, and such that 
\begin{itemize}
\item[(i)]  the critical points of $f_t$ are at 0 and $\infty$ for all $t\in (0,1]$, 
\item[(ii)]  the critical point at 0 is in a cycle $0 \mapsto -1 \mapsto 0$ 
for all $t$,  and
\item[(iii)]  there is an attracting fixed point of multiplier $\alpha(t)\to -1$
as $t\to 0$.  
\end{itemize}
Then this family is contained in the hyperbolic component of 
$f_1$.  Recall that the Julia set of $f_1$ is the basilica.  
The family $f_t$ can be expressed as 
  $$f_t(z) = \frac{z^2 - 1}{c(t) z^2 + 1},$$
for a function $c$ with $c(1) = 0$.  It can be computed directly
that the triple of fixed point multipliers of $f_t$ tends to 
$\{\infty, -1, -1\}$ as $c(t)$ descends from 0 to $-1$.  In the limit, the second iterate
of $f_t$ converges to $-2z^2/(z^2 + 1)$ locally uniformly on $\Chat -
\{-1, 1\}$, which is conjugate to $G_{-1}(z) = z - 1 + 1/z$ and 
to the polynomial $z^2 + 1/4$ with
a parabolic fixed point.  It follows that the $\tau^2$-value (as defined
before the statement of Theorem \ref{tau2}) for this family is equal to 1.  

The Julia sets of $f_t$ appear to 
converge (in the Hausdorff topology) to the 
cauliflower Julia set of $z^2 + 1/4$.  The divergence of this family $[f_t]$ in $M_2$
is an illustration of the obstruction to mating the polynomial $f_1$ with 
itself (the $1/2$-limb of the Mandelbrot set is its own conjugate).


\bigskip\bigskip

\begin{thebibliography}{MFK}

\bibitem[De]{D:measures}
L.~DeMarco.
\newblock {Iteration at the boundary of the space of rational maps}.
\newblock {\em Submitted to Duke Math. Journal, \em{March, 2004}}.

\bibitem[DE]{Douady:Earle}
A.~Douady and C.~Earle.
\newblock {Conformally natural extension of homeomorphisms of the circle}.
\newblock {\em Acta Math.} {\bf 157}(1986), 23--48.

\bibitem[Ep]{Epstein:bounded}
A.~Epstein.
\newblock {Bounded hyperbolic components of quadratic rational maps}.
\newblock {\em Ergodic Theory Dynam. Systems} {\bf 20}(2000), 727--748.

\bibitem[FLM]{FLM}
A.~Freire, A.~Lopes, and R.~Ma\~{n}{\'e}.
\newblock {An invariant measure for rational maps}.
\newblock {\em Bol. Soc. Brasil. Mat.} {\bf 14}(1983), 45--62.

\bibitem[Ly]{Lyubich:entropy}
M.~Lyubich.
\newblock {Entropy properties of rational endomorphisms of the {R}iemann
  sphere}.
\newblock {\em Ergodic Theory Dynamical Systems} {\bf 3}(1983), 351--385.

\bibitem[Ma1]{Mane:unique}
R.~Ma\~{n}{\'e}.
\newblock {On the uniqueness of the maximizing measure for rational maps}.
\newblock {\em Bol. Soc. Brasil. Mat.} {\bf 14}(1983), 27--43.

\bibitem[Ma2]{Mane}
R.~Ma\~{n}{\'e}.
\newblock {The {H}ausdorff dimension of invariant probabilities of rational
  maps}.
\newblock In {\em Dynamical Systems, Valparaiso 1986}, pages 86--117. Springer,
  Berlin, 1988.

\bibitem[Mi]{Milnor:quad}
J.~Milnor.
\newblock {Geometry and dynamics of quadratic rational maps}.
\newblock {\em Experiment. Math.} {\bf 2}(1993), 37--83.
\newblock With an appendix by the author and Lei Tan.

\bibitem[MFK]{Mumford:GIT}
D.~Mumford, J.~Fogarty, and F.~Kirwan.
\newblock {\em Geometric invariant theory}, volume~34 of {\em Ergebnisse der
  Mathematik und ihrer Grenzgebiete (2) [Results in Mathematics and Related
  Areas (2)]}.
\newblock Springer-Verlag, Berlin, third edition, 1994.

\bibitem[Si]{Silverman}
J.~H. Silverman.
\newblock {The space of rational maps on $\bf {P}\sp 1$}.
\newblock {\em Duke Math. J.} {\bf 94}(1998), 41--77.

\end{thebibliography}

\def\cprime{$'$}

\end{document}